\documentclass[12pt,centertags,oneside]{article}
\usepackage{amsmath,amstext,amsthm,amscd,typearea}
\usepackage{amssymb}
\usepackage{a4wide}
\usepackage[mathscr]{eucal}
\usepackage{mathrsfs}
\usepackage{typearea}
\usepackage{charter}
\usepackage{pdfsync}
\usepackage{url}

\usepackage{xcolor}

\usepackage[a4paper,width=16.2cm,top=3cm,bottom=3cm]{geometry}

\numberwithin{equation}{section}

\newtheorem{theorem}{Theorem}[section]
\newtheorem{definition}[theorem]{Definition}
\newtheorem{proposition}[theorem]{Proposition}
\newtheorem{corollary}[theorem]{Corollary}
\newtheorem{lemma}[theorem]{Lemma}
\newtheorem{remark}[theorem]{Remark}

\newtheorem{question}[theorem]{Question}

\newcommand{\cali}[1]{\mathscr{#1}}

\newcommand{\supp}{{\rm Supp}}

\newcommand{\ddc}{dd^c}
\newcommand{\dc}{d^c}

\newcommand{\codim}{{\rm codim\ \!}}

\newcommand{\capK}{\text{cap}}

\newcommand{\C}{\mathbb{C}}

\newcommand{\N}{\mathbb{N}}

\newcommand{\R}{\mathbb{R}}

\title{\bf  Relative non-pluripolar product of currents}
\providecommand{\keywords}[1]{\textbf{\textit{Keywords:}} #1}
\providecommand{\subject}[1]{\textbf{\textit{Mathematics Subject Classification 2010:}} #1}

\author{Duc-Viet Vu}

\newcommand{\Addresses}{{
		\bigskip
		\footnotesize
		\textsc{Duc-Viet Vu, University of Cologne, Division of Mathematics, Department of Mathematics and Computer Science, Weyertal 86-90, 50931, K\"oln,  Germany  \& Thang Long Institute of Mathematics and Applied Sciences, Hanoi, Vietnam}
		\noindent
		\par\nopagebreak
		\noindent
		\textit{E-mail address}: \texttt{vuduc@math.uni-koeln.de}	
}}

\date{\today}
\begin{document}
\maketitle
\begin{abstract} Let $X$ be a compact K\"ahler manifold. Let $T_1, \ldots, T_m$ be closed positive currents of bi-degree $(1,1)$ on $X$ and $T$ an arbitrary closed positive current on $X$. We introduce the \emph{non-pluripolar product relative to $T$} of $T_1, \ldots, T_m$.  We recover the well-known non-pluripolar product of $T_1, \ldots, T_m$ when $T$ is the current of integration along $X$. Our main results are a monotonicity property of relative  non-pluripolar products, a necessary condition for currents to be of relative full mass intersection in terms of Lelong numbers, and the convexity of weighted classes of currents of  relative full mass intersection. The former two results are new even when $T$ is the current of integration along $X$.  
\end{abstract}
\noindent
, {quasi-continuity}, {weighted class},  {Monge-Amp\`ere operator}, {full mass intersection}.
\\

\noindent
\subject{32U40}, {32H50}, {37F05}.

\tableofcontents


\section{Introduction}

The problem of defining the intersection of  closed positive currents on a complex manifold   is a central question in the pluripotential theory. This question is already important for currents of bi-degree $(1,1)$ with  deep applications to complex geometry as well as complex dynamics; see \cite{Bedford_Taylor_82,Demailly_ag,Kolodziej05,Fornaess_Sibony,Dinh_Sibony_density} for an introduction to the subject. 

Let $X$ be an arbitrary complex manifold. Let $T_1, \ldots,T_m$ be closed positive currents of bi-degree $(1,1)$ on $X$.   In \cite{BT_fine_87,GZ-weighted,BEGZ}, the  non-pluripolar product $\langle T_1 \wedge \cdots \wedge T_m\rangle$ of $T_1, \ldots, T_m$ was defined  by using the quasi-continuity of plurisubharmonic (psh for short) functions with respect to capacity. Since then, this notion has played an important role  in complex geometry. We refer to \cite{Darvas_book,Lu-Darvas-DiNezza-mono,Lu-Darvas-DiNezza-singularitytype,Lu-Darvas-DiNezza-logconcave} for recent developments.  

In this paper,  our first goal is to generalize the notion of non-pluripolar product to a more general natural setting. Give a closed positive current $T$ on $X$,   we introduce \emph{the non-pluripolar product relative to $T$} of $T_1, \ldots, T_m$ which we denote by  $\langle T_1 \wedge \cdots \wedge T_m \dot{\wedge} T \rangle $.  We recover the above-mentioned non-pluripolar product when $T$ is the current of integration along $X$, see also Remark \ref{re-currentinterV}.  When well-defined, $\langle T_1 \wedge \cdots \wedge T_m \dot{\wedge} T \rangle $  is a closed positive current on $X$. As in \cite{BEGZ}, the proof of the last fact follows from ideas in \cite{Sibony_duke}. The product $\langle T_1 \wedge \cdots \wedge T_m \dot{\wedge} T \rangle $ is symmetric, homogeneous and sub-additive in $T_1, \ldots, T_m$ (Proposition \ref{pro-sublinearnonpluripolar}). One can consider the last product as a sort of intersection of $T_1, \ldots, T_m,T$.

The main ingredient in  the construction of relative non-pluripolar products is  a  \emph{(uniform) quasi-continuity} property of \emph{bounded} psh functions with respect to the capacity associated to closed positive currents (Theorem \ref{th-capconvergedecreasing}). This quasi-continuity property is stronger than the usual one. We also need to establish some convergence properties of mixed Monge-Amp\`ere operators which are of independent interest.

Consider now the case where $X$ is a compact K\"ahler manifold. 
As one can expect, the relative non-pluripolar products are always well-defined in this setting.  For any closed positive current $R$ in $X$, we denote by $\{R\}$ the cohomology class of $R$. For cohomology classes $\alpha,\beta$ in $H^{p,p}(X,\R)$, we write $\alpha \le \beta$ if $(\beta-\alpha)$ can be represented by a closed positive current. The following result gives a monotonicity property of relative non-pluripolar products.

\begin{theorem} \label{th-main1} (Theorem \ref{th-monoticity}) Let $X$ be  a compact K\"ahler manifold and $T_1,\ldots, T_m,T$ closed positive currents on $X$ such that $T_j$ is of bi-degree $(1,1)$ for $1 \le j \le m$.  Let $T'_j$ be closed positive $(1,1)$-current in the cohomology class of $T_j$ on $X$ such that $T'_j$ is less singular than $T_j$ for $1 \le j \le m$. Then we have 
$$\{\langle T_1 \wedge \cdots \wedge T_m \dot{\wedge} T \rangle \} \le \{\langle T'_1 \wedge \cdots \wedge T'_m \dot{\wedge} T \rangle\}.$$
\end{theorem}

Note that even when $T$ is the current of integration along $X$, this result is new.  In the last case,  the above result was conjectured in \cite{BEGZ} and  proved there for $T_j$ of potentials locally bounded outside a closed complete pluripolar set. The last condition was relaxed in \cite{Lu-Darvas-DiNezza-mono,WittNystrom-mono} but it was required there that  $m=\dim X$ (always for $T$ to be the current of integration along $X$). The proofs of these results presented there don't extend to our setting.  A key ingredient in our proof is Theorem \ref{th-increasingsequenceMa} (and Remark \ref{re-th-increasingsequenceMa}) giving a generalization of well-known convergence properties of Monge-Amp\`ere operators. To prove  Theorem \ref{th-increasingsequenceMa}, we will need the strong quasi-continuity of \emph{bounded} psh functions mentioned above.  When $T$ is of bi-degree $(1,1)$, we actually have a stronger monotonicity property, see Remark \ref{re-Tofbidegre11mono}. 


Theorem \ref{th-main1} allows us to define the notion of \emph{full mass intersection relative to $T$} for currents $T_1, \ldots, T_m$ as in \cite{BEGZ}, see Definition \ref{def-fullmassinter}.  Our next goal is to study currents of  relative full mass intersection.  We will prove that the relative non-pluripolar products of currents with full mass intersection are continuous under decreasing or increasing sequences (Theorem \ref{th-convergencMAgenernonpluri}).   

Now we concentrate on  the case where the cohomology classes of $T_1, \ldots, T_m$  are K\"ahler. In this case, $T_1, \ldots, T_m$ are of full mass intersection relative to $T$ if and only if $\{ \langle \wedge_{j=1}^m T_j \dot{\wedge} T\rangle\}$ is equal to $\{\wedge_{j=1}^m \theta_j \wedge T\}$, where $\theta_j$ is a closed smooth $(1,1)$-form cohomologous to $T_j$ for $1 \le j \le m$.  
 Our next main result explains why having positive Lelong numbers is an obstruction for  currents to be of  full mass intersection. 

\begin{theorem} \label{th-main3} (Theorem \ref{the-lelongobstructiontongmpnhonhonn})   Let $X$ be a compact K\"ahler  manifold. Let $T_1, \ldots, T_m$ be closed positive $(1,1)$-currents on $X$ such that the cohomology class of $T_j$ is K\"ahler for $1 \le j \le m$  and $T$ a closed positive $(p,p)$-current on $X$ with $p+m \le n$.   Let $V$ be an irreducible analytic subset in $X$ such that the generic Lelong numbers of $T_1, \ldots, T_m,T$ along $V$ are strictly positive. Assume $T_1, \ldots, T_m$ are of full mass intersection relative to $T$. Then, we have 
$\dim V <n-p-m$. 
\end{theorem}

Recall that the generic Lelong number of $T$ along $V$ is the smallest value among the Lelong numbers of $T$ at points in $V$.  When $p+m=n$, the above theorem says that $V$ is empty.  To get motivated about Theorem \ref{th-main3},  let  us consider the case where $p+m=n$ and  the current $T_j$ has a potential which is locally bounded outside a point $x_0$ in $X$. In this case,  $T_1 \wedge \cdots \wedge T_m \wedge T$ is  well-defined in the sense given in \cite{Bedford_Taylor_82,Demailly_ag,Fornaess_Sibony} and we can see that $\langle T_1 \wedge \cdots \wedge T_m \dot{\wedge} T \rangle = T_1 \wedge \cdots \wedge T_m \wedge T$ on $X \backslash \{x_0\}$. The latter current has strictly positive mass on $x_0$ by \cite[Corollary 7.9]{Demailly_ag}. Hence $T_1, \ldots, T_m$ cannot be of full mass intersection relative to $T$. By considering $T_j$ to be  suitable currents with analytic singularities (for example when $X$ is projective) and $T$ is the current of integration along $X$, one can see that  the conclusion of Theorem \ref{th-main3} is optimal.  

When $m=n$ and $T_1= \cdots = T_m$, Theorem \ref{th-main3} was proved  in \cite{GZ-weighted}. Their proof uses comparison principle and hence doesn't apply to our setting because $T_1, \ldots, T_m$ are different in general. 
When the class of $T_j$ is not K\"ahler, the above theorem no longer holds because currents with minimal singularities in a big and non-nef cohomology class always have positive Lelong numbers at some points in $X$ (see \cite{Boucksom-these}).  The proof of Theorem \ref{th-main3}, which is Theorem \ref{the-lelongobstructiontongmpnhonhonn} below,  uses Theorem \ref{th-main1}. 

Let $\mathcal{W}^-$ be the set of convex increasing functions $\chi: \R \to \R$ with $\chi(-\infty)= -\infty$. A function in $\mathcal{W}^{-}$ is called a \emph{(convex) weight}.  We will define the notion of having \emph{full mass intersection relative to $T$ with weight $\chi$} and the \emph{weighted class} of  these currents, see Definition \ref{def-weigtedclasschi}.  The case where $T$ is the current of integration along $X$ and $T_1, \ldots, T_m$ are all equal  was studied in \cite{GZ-weighted,BEGZ}.   Here is our main result in this part.

\begin{theorem} \label{th-monoconvexweight} Let $X, T_1,\ldots, T_m,T$ be as in Theorem \ref{th-main3}. Let $\chi \in \mathcal{W}^-$. The following three assertions hold:

$(i)$ If  $T_j, \ldots, T_j$  ($m$ times $T_j$) are of full mass intersection relative to $T$  for $1 \le j \le m$, then  $T_1, \ldots, T_m$ are also of full mass intersection relative to $T$.

$(ii)$ Let  $T'_j $ be  a closed positive $(1,1)$-current whose cohomology class is K\"ahler  for $1 \le j \le m$. Assume that $T'_j$ is less singular than $T_j$ and $T_1, \ldots, T_m$ are of full mass intersection relative to $T$ with weight $\chi$. Then $T'_1, \ldots, T'_m$ are also of full mass intersection relative to $T$ with weight $\chi$.

$(iii)$ If $T_j, \ldots, T_j$  ($m$ times $T_j$) are of full mass intersection relative to $T$  with weight $\chi$   for $1 \le j \le m$, then  $T_1, \ldots, T_m$ are also of full mass intersection relative to $T$ with weight $\chi$.
\end{theorem}

Theorem \ref{th-monoconvexweight} is  a combination of Theorems \ref{the-mathcalEnalpha},  \ref{the-monotonicityforweigtedclass} and \ref{the-mathcalEnalphaweighted} below. 
Property $(i)-(ii)$ of Theorem \ref{th-monoconvexweight} was proved in \cite{GZ-weighted} in the case where $T$ is the current of integration along $X$ and $T_1, \ldots, T_m$ are  equal and $T'_1, \ldots, T'_m$ are equal.  When $T$ is the current of integration along $X$, Property $(iii)$ of  Theorem \ref{th-monoconvexweight} was proved in \cite{GZ-weighted} for $\chi(t)= -(-t)^r$ with $0< r \le 1$ and it was proved in \cite{Darvas_book} for every $\chi \in \mathcal{W}^-$ under an extra condition that $m=n$. The proof given in \cite{Darvas_book}  used Monge-Amp\`ere equations, hence, doesn't extend to our setting. 

We have some comments about the proof of Theorem \ref{th-monoconvexweight}. Firstly,  arguments similar to those in \cite{GZ-weighted} are sufficient to obtain $(ii)$ of Theorem \ref{th-monoconvexweight}. However, as pointed out above, to prove $(iii)$, we need new arguments.  We prove $(i), (iii)$ using the same approach.  It is possible that we can use the comparison principle to get $(i)$ as in the case where $T$ is the current of integration along $X$. But we choose to invoke a more flexible idea which is also applicable to get $(iii)$. The idea is to combine a monotonicity property  and generalizations of results in \cite{BT_fine_87} about plurifine topology properties of Monge-Amp\`ere operators. Here, the monotonicity mentioned in the last sentence is Theorem \ref{th-main1} in case of proving $(i)$ and $(ii)$ in case of proving $(iii)$.

The paper is organized as follows. In Section \ref{sec-quasicontinuity}, we establish the uniform strong quasi-continuity of bounded psh functions and derive from it some consequences concerning the convergence of Monge-Amp\`ere operators.  In Section \ref{sec-def-gene-nonpluri}, we define the notion of relative non-pluripolar products and prove its basic properties. In Section \ref{sec-mono}, we prove the monotonicity of relative non-pluripolar products and explain Lelong number obstruction for currents having full mass intersection.  Section \ref{sec-weightedclass} is devoted to the weighted class of currents with relative full mass intersection. \\

\noindent
\textbf{Acknowledgments.}  The author would like to thank Hoang Chinh Lu  for kindly answering him numerous questions about non-pluripolar products. He also thanks Tam\'as Darvas, Lucas Kaufmann and Tuyen Trung Truong for their comments. This research  is supported by a postdoctoral fellowship of the Alexander von Humboldt Foundation. \\

\noindent
\textbf{Notation and convention.} 
For a closed positive current $T$ on a compact complex manifold $X$, we denote by $\{T\}$ the cohomology class of $T$ and $\nu(T,x)$ denotes the Lelong number of $T$ at $x$.  Recall that the wedge product on closed  smooth forms induces the cup-product on their de Rham cohomology  classes. We use the same notation $`` \wedge "$ to denote these two products.

Recall $\dc=  \frac{i}{2\pi} (\overline \partial - \partial)$ and  $\ddc = \frac{i}{\pi} \partial \overline \partial$.  We use $\gtrsim, \lesssim$ to denote $\ge, \le$ modulo some multiplicative constant independent of parameters in question. 
For a  Borel  set $A \subset X$, we denote by $\bold{1}_A$ the characteristic function of $A$, that means $\bold{1}_A$ is equal to $1$ on $A$ and $0$ elsewhere. For every current $R$ of order $0$ on $X$, we  denote by $\|R\|_A$ the mass of $R$ on $A$.


\section{Quasi-continuity for bounded psh functions} \label{sec-quasicontinuity}

In this section, we prove a quasi-continuity for bounded plurisubharmonic (psh for short) functions which is stronger than the one for general psh functions. This property is the key to define our generalization of non-pluripolar products. 

Let $U$ be an open subset of $\C^n$.   
Let $K$ be a Borel subset of  $U$.  The  \emph{capacity} $\capK(K,U)$ of $K$ in $U$, which was introduced in \cite{Bedford_Taylor_82}, is given by 
$$\capK(K,U):= \sup \bigg\{ \int_K (\ddc u)^n:  u \text{ is psh  on } U \text{ and } 0 \le u \le 1\bigg\}.$$
For a  closed positive current $T$ of bi-dimension $(m,m)$ on $U$ ($0 \le m \le n$),  we define
$$\capK_T(K,U):= \sup \bigg\{ \int_K (\ddc u)^m \wedge T:  u \text{ is psh  on } U \text{ and } 0 \le u \le 1\bigg\}.$$
We say that a sequence of functions $(u_j)_{j \in \N}$ \emph{converges to $u$ with respect to the capacity $\capK_T$ (relative in $U$)} if for any constant $\epsilon>0$ and $K \Subset U$, we have $\capK_T\big(\{|u_j - u| \ge \epsilon\}\cap K, U\big) \to 0$ as $t \to \infty$.  The last notion was introduced in \cite{Kolodziej05,Xing-continuity}.  We say that a subset $A$ in $U$ is \emph{locally complete pluripolar} set if  locally $A= \{\psi= -\infty\}$ for some psh function $\psi$. We begin with the following lemma which is probably well-known. 

\begin{lemma} \label{le-capTK} Let $A$ be a locally complete pluripolar set in $U$. Let $T$ be closed positive current $T$ of bi-dimension $(m,m)$ on $U$. Assume that  $T$ has no mass on $A$. Then,  we have $\capK_T(A,U)=0$. 
\end{lemma}

\proof The proof is standard. We present the details for readers' convenience. Since the problem is of local nature,  we can assume that there is a negative psh function $\psi$ on $U$ such that $A= \{\psi= -\infty\}$.  Let $u_1,\ldots, u_m$ be bounded psh functions on $U$ such that $0 \le u_j \le 1$ for $1 \le j \le m$.  
Let $\omega$ is the standard K\"ahler form on $\C^n$. Let $k \in \N$ and $\psi_k:= k^{-1} \max \{\psi, -k\}$. We have $-1 \le \psi_k \le 0$. 
Let  $\chi$ be a nonegative smooth function  with compact support in $U$. Let $0 \le l \le m$ be an integer.  Put 
 $$I_k:=  \int_{U} \chi \psi_k \ddc u_1 \wedge \cdots \wedge \ddc u_l \wedge \omega^{m-l} \wedge T.$$
  Since $\psi_k =-1$ on $\{\psi <-k\}$, in order to prove the desired assertion, it is enough to show that for every $0 \le l \le q$,  we have
\begin{align}\label{eq-capKATlw-CAPtk2}
 I_k \to 0
\end{align}
as $k \to\infty$ uniformly in $u_1,\ldots, u_l$. We will prove  (\ref{eq-capKATlw-CAPtk2}) by induction on $l$. Firstly, (\ref{eq-capKATlw-CAPtk2}) is trivial if $l=0$ because $T$ has no mass on $A$. Assume that it holds for $(l-1)$. We prove it for $l$.  Put 
$$R:= \ddc u_2 \wedge \cdots \wedge \ddc u_l \wedge \omega^{m-l} \wedge T.$$ 
By integration by parts, we have
$$I_k= \int_U u_1  \chi \ddc \psi_k \wedge R + \int_U u_1 \psi_k \ddc \chi \wedge R+2 \int_U u_1 d \psi_k \wedge \dc \chi \wedge  R.$$
Denote by $I_{k,1}, I_{k,2},  I_{k,3}$ the first, second and third term respectively in the right-hand side of the last equality. Since $u_1$ is bounded by $1$, by integration by parts, we get 
$$|I_{k,1}| \le C \int_{\supp \chi} -\psi_k R \wedge \omega, \quad |I_{k,2}| \le C \int_{\supp \chi}- \psi_k R \wedge \omega,$$
for some constant $C$ depending only on $\chi$. By induction hypothesis, we have 
$$\lim_{k \to \infty} \int_{\supp \chi} \psi_k R \wedge \omega= 0.$$
 Thus $\lim_{k \to \infty} I_{k,j} =0$ for $j=1,2$. To treat $I_{k,3}$, we use the Cauchy-Schwarz inequality to get
 $$|I_{k,3}| \le \bigg(\int_{\supp \chi} d\psi_k  \wedge \dc \psi_k \wedge R\bigg)^{1/2} \ \lesssim \bigg(\int_{U_1} -\psi_k  R \wedge \omega\bigg)^{1/2},$$
because $d \psi_k \wedge \dc \psi_k= \ddc \psi_k^2 - \psi_k \ddc \psi_k$,  where $U_1$ is a relatively compact open subset of $U$ containing the support of $\chi$. By induction hypothesis, $\lim_{k \to \infty}\int_{U_1} \psi_k  R \wedge \omega= 0$. So $\lim_{k\to \infty} I_{k,3}=0$. In conclusion, (\ref{eq-capKATlw-CAPtk2}) follows. This finishes the proof.
\endproof 

We now give a definition which will be crucial later. Let $(T_k)_k$ be a sequence of closed positive currents of bi-degree $(m,m)$ on $U$. We say that $(T_k)_k$ satisfies \emph{Condition $(*)$} if $(T_k)_k$ is of uniformly bounded mass on compact subsets of $U$, and for every open set $U' \subset U$ and every bounded psh function $u$ on $U'$ and every sequence $(u_k)_k$ of psh functions on $U'$ decreasing to $u$, we have 
\begin{align}\label{inere-chuyuniformquasicontinuiTltongquadef}
\lim_{k \to \infty}(u_k- u) (\ddc u)^{m} \wedge T_{k}= 0
\end{align}
An obvious example for sequences satisfying Condition $(*)$ is constant sequences: $T_k= T$ for every $k$. We can also take $(T_k)_k$ to be a sequence of suitable Monge-Amp\`ere operators with decreasing potentials, see \cite{Bedford_Taylor_82,Demailly_ag,Fornaess_Sibony}. Concerning Condition $(*)$, we will use mostly the example of constant sequences and the one provided by the following result. 

\begin{theorem} \label{th-decreasing-unbounded}  Let $S$ be a closed positive current on $U$. Let $v$ be a psh function on $U$ such that $v$ is locally integrable with respect to the trace measure of $S$ and $(v_k)_k$ a sequence of psh functions on $U$ such that $v_k \to v$ in $L^1_{loc}$  as $k \to \infty$ and $v_k \ge v$ for every $k$. Let $T:= \ddc v\wedge S$ and $T_k:= \ddc v_k \wedge S$.  Let $u_j$ be a bounded psh function on $U$ for $1 \le j \le m$. Let $(u_{jk})_{k \in \N}$ be  a sequence of  uniformly  bounded psh functions such that  $u_{jk} \to u_j$ in $L^1_{loc}$ as $k \to \infty$ and $u_{jk} \ge u_j$ for every $j,k$. Then we have 
\begin{align}\label{limit-decreaseingpointwise}
u_{1k} \ddc u_{2k} \wedge \cdots \wedge \ddc u_{mk} \wedge T_k  \to u_1 \ddc u_{2} \wedge \cdots \wedge \ddc u_{m} \wedge T
\end{align}
as $k \to \infty$. In particular, the sequence $(T_k)_k$ satisfies Condition $(*)$. 
\end{theorem}    

\proof By Hartog's lemma, $v_k, u_{jk}$ are uniformly bounded from above in $k$ on compact subsets of $U$ for every $j$.  Since the problem is local, as usual, we can assume that $U$ is relatively compact open set with smooth boundary in $\C^n$, every psh function in questions is defined on an open neighborhood of $\overline U$, $v_k, v \le 0$ on $U$ for every $k$ and  $u_{jk}, u_j$ are all equal to a smooth psh function $\psi$ outside some fixed compact subset of $U$ such that $\psi =0$ on $\partial U$. 
We claim that
\begin{align} \label{limit-decreaseingpointwisevodau}
Q_k:=v_k \ddc u_{1k} \wedge \cdots \wedge \ddc u_{m k}\wedge S \to v \ddc u_{1} \wedge \cdots \wedge \ddc u_{m} \wedge S
\end{align}
 as $k \to\infty$. In particular, this implies that $v$ is locally integrable with respect to $\ddc u_{1} \wedge \cdots \wedge \ddc u_{m}\wedge S$.  We will prove (\ref{limit-decreaseingpointwise}) and (\ref{limit-decreaseingpointwisevodau}) simultaneously  by induction on $m$.  When $m=0$, we have nothing to prove. Assume that (\ref{limit-decreaseingpointwise}) and (\ref{limit-decreaseingpointwisevodau}) hold for $(m-1)$ in place of $m$.  Let 
$$R_{j,k}:= \ddc u_{jk} \wedge \cdots \wedge \ddc u_{mk} \wedge T_k$$
for $1 \le j \le m$.  By induction hypothesis, we have 
$$R_{j,k} \to R_j:=  \ddc u_{j} \wedge \cdots \wedge \ddc u_{m} \wedge T$$
for $j \ge 2$. Since $u_{1k}$ is uniformly bounded on  $U$, the family $u_{1k} R_{2,k}$ is of uniformly bounded mass. Let $R_\infty$ be a limit current of the last family. Without loss of generality, we can assume $R_\infty= \lim_{k \to \infty} u_{1k} R_{2,k}$ and $S$ is of bi-degree $(n-m,n-m)$. By standard arguments, we have  $R_\infty \le u_1 R_2$ (see \cite[Proposition 3.2]{Fornaess_Sibony}). Thus, in order to have $R_\infty = u_1 R_2$, we just need to check that 
\begin{align}\label{eq-massth-decreasing-unbounded}
\int_U  R_\infty \ge \int_U u_1 R_2
\end{align}
 (both sides are finite because of the assumption we made at the beginning of the proof). Since $\psi= 0$ on $\partial U$ and $u_{1k}=\psi$ on outside a compact of  $U$, we have 
\begin{align}\label{limit-trenUpshu1kR2k}
\int_U u_{1k} R_{2,k} \to \int_U R_\infty , \quad \int_U \psi R_{2,k}  \to \int_U \psi R_2.
\end{align}
Let $u_{jk}^\epsilon, \psi^\epsilon$ be standard regularisations of $u_{jk}, \psi$ respectively.  Since $u_{jk}= \psi$ outside some  compact of $U$, we have $u_{jk}^\epsilon = \psi^\epsilon$ outside some compact $K$ of $U$, for $\epsilon$ small enough and $K$ independent of $j,k,\epsilon$. Consequently, $u_{jk}^\epsilon - \psi^\epsilon$ is  supported in $K \Subset U$. Note that since $\psi$ is smooth, $\psi^\epsilon \to  \psi$ in $\cali{C}^\infty$- topology.  By  integration by parts and the fact that $u_{jk} \ge u_j$ for $j=1,2$, we have 
\begin{align*}
\int_U (u_{1} -\psi) R_{2} &\le  \lim_{\epsilon \to 0} \int_U (u^\epsilon_{1k} -\psi^\epsilon) R_{2} = \lim_{\epsilon \to 0}\int_U u_{2}  \ddc (u^\epsilon_{1k} -\psi^\epsilon) R_{3}\\
&\le \lim_{\epsilon \to 0}\int_U u^\epsilon_{2k}  \ddc (u^\epsilon_{1k} -\psi^\epsilon) R_{3}+ \lim_{\epsilon \to 0}\int_U (u^\epsilon_{2k}- u_{2}) \ddc \psi^\epsilon \wedge  R_{3}\\
& =\lim_{\epsilon \to 0}\int_U   (u^\epsilon_{1k} -\psi^\epsilon)\ddc u^\epsilon_{2k} \wedge   R_{3} + o_{k \to \infty}(1)  
\end{align*}
by induction hypothesis for $(m-1)$ of (\ref{limit-decreaseingpointwise}) and the fact that $\| \ddc \psi_\epsilon - \ddc \psi\|_{\cali{C}^0}= O(\epsilon)$. Put $R'^\epsilon_{2,k}:= \ddc u^\epsilon_{2k} \wedge \cdots \wedge \ddc u^\epsilon_{mk}$.   Repeating the above arguments for every $u_{jk}$ ($j \ge 2$) and $v,v_k$  gives
\begin{align*}
\int_U (u_{1} -\psi) R_{2} &\le \lim_{\epsilon \to 0}\int_U   (u^\epsilon_{1k} -\psi^\epsilon) R'^\epsilon_{2,k} \wedge \ddc v \wedge S  \le  \lim_{\epsilon \to 0}\int_U   v \ddc (u^\epsilon_{1k} -\psi^\epsilon)\wedge R'^\epsilon_{2,k} \wedge S  \\
&\le \lim_{\epsilon \to 0}\int_U   v_k \ddc (u^\epsilon_{1k} -\psi^\epsilon) \wedge R'^\epsilon_{2,k} \wedge S+   \lim_{\epsilon \to 0}\int_U   (v_k- v) \ddc \psi^\epsilon \wedge R'^\epsilon_{2,k} \wedge S\\
&= \lim_{\epsilon \to 0}\int_U  (u^\epsilon_{1k} -\psi^\epsilon) \wedge R'^\epsilon_{2,k} \wedge \ddc v_k \wedge S+ o_{k \to \infty}(1)\\
&=\int_U  (u_{1k} -\psi) \wedge R_{2,k}+ o_{k \to \infty}(1) 
\end{align*}
by  (\ref{limit-decreaseingpointwisevodau}) for $(m-1)$ and the usual convergence of Monge-Amp\`ere operators. Letting $k \to \infty$ in the last inequality and using (\ref{limit-trenUpshu1kR2k}) give (\ref{eq-massth-decreasing-unbounded}). Hence (\ref{limit-decreaseingpointwise}) for $m$ follows.

It remains to prove   (\ref{limit-decreaseingpointwisevodau}) for $m$. Put $R'_{2,k}:=  \ddc u_{2k} \wedge \cdots \wedge \ddc u_{m k}$ and $R'_2:=  \ddc u_{2} \wedge \cdots \wedge \ddc u_{m}$. We check that $Q_k$ is of uniformly bounded mass.  Decompose
\begin{align*}
Q_k =v_k \ddc (u_{1k}- \psi) \wedge R'_{2,k}\wedge S+ v_k \ddc \psi \wedge R'_{2,k}\wedge S. 
\end{align*}
The second term converges to $v \ddc \psi \wedge R'_{2}\wedge S$ as $k \to \infty$ by induction hypothesis for $(m-1)$. Denote by $Q_{k,1}$ the first term. Let $v^\epsilon_k$ be standard regularizations of $v_k$. By integration by parts, we have 
\begin{align*}
\int_U Q_{k,1}^\epsilon &:=  \int_U v^\epsilon_k \ddc (u_{1k}- \psi) \wedge R'_{2,k}\wedge S\\
&= \int_U (u_{1k}- \psi)  \ddc v^\epsilon_k \wedge R'_{2,k} \wedge S=(u_{1k}- \psi) R'_{2,k}\wedge   \ddc v^\epsilon_k \wedge S
\end{align*}
which converges to $\int_U (u_{1k}- \psi) R'_{2,k} \wedge \ddc v_k\wedge S$ as $\epsilon \to 0$  by (\ref{limit-decreaseingpointwise}) for $m$. Thus, 
$$\int_U Q_{k,1}=  \int_U (u_{1k}- \psi) R'_{2,k} \wedge \ddc v_k\wedge S.$$
This combined with (\ref{limit-decreaseingpointwise}) for $m$ again implies that  $\int_U Q_{k,1} \to \int_U (u_1 - \psi) R'_2 \wedge \ddc v \wedge S$ as $k \to \infty$.  The last limit is equal to $\int_U v \ddc (u_1 - \psi) \wedge R'_2$ by integration by parts which can be performed thanks to (\ref{limit-decreaseingpointwise}) for $m$.  Thus, we have proved that $Q_k$ is of uniformly bounded mass and 
$$\int_U Q_k \to \int_U v R_1$$
as $k \to \infty$. This combined with the fact that $v R_1 \ge Q_\infty$ for every limit current $Q_\infty$ of the family $(Q_k)_k$ gives  the desired assertion  (\ref{limit-decreaseingpointwisevodau}) for $m$.    This finishes the proof.
\endproof

\begin{remark} \label{re-symmetricMA} We can apply Theorem \ref{th-decreasing-unbounded} to the case where $S$ is a constant function and $v$ is an arbitrary psh function.  By the above proof, we can check the following observations: 

$(i)$  if $S$ is a closed positive current on $U$, $u_1, \ldots, u_m$ are psh functions on $U$ which are locally integrable with respect to $S$ such that $u_j$ is locally bounded for every $1 \le j \le m$  except possibly for  one index, then the current $\ddc u_1 \wedge \cdots \wedge \ddc u_m\wedge S$, which is defined inductively as usual,  is symmetric with respect to $u_1, \ldots, u_m$ and satisfies the convergence under decreasing sequences,

$(ii)$  let $u_0$ be another psh function locally integrable with respect to $S$ such that  $u_0$ is locally bounded if there is an index $1 \le j \le m$ so that $u_j$ is not locally bounded. Then $u_0 \ddc u_1 \wedge \cdots \wedge \ddc u_m \wedge S$ is convergent under decreasing sequences and for every compact $K$ in $U$, if we have $0\le u_1,\ldots, u_m\le 1$, then  
$$\|u_0 \ddc u_1 \wedge \cdots \wedge \ddc u_m \wedge S\|_K \le C \|u_0 S\|_U$$
for some constant $C$ independent of $u_0, \ldots, u_m,S$. 
\end{remark}

The following result explains the reason for the use of Condition $(*)$.

\begin{theorem}\label{th-capconvergedecreasing}  Let $(T_l)_l$ be a sequence of closed positive currents satisfying Condition $(*)$.  Let $u$ be a \emph{bounded}  psh function on $U$ and  $(u_k)_k$ a sequence of psh functions on $U$ decreasing to $u$.    Then for every constant $\epsilon >0$ and every compact $K$ in $U$, we have $\capK_{T_l}(\{|u_k -u| \ge \epsilon\}\cap K) \to 0$ as $k \to \infty$ uniformly in $l$. In particular,  for every constant $\epsilon>0$,  there exists an open subset $U'$ of $U$ such that  $\capK_{T_l}(U',U)<\epsilon$ for every $l$ and the restriction of  $u$ to $U \backslash U'$ is continuous. 
\end{theorem}

Consider the case where $T_l= T$ for every $l$. Then, the above theorem give a quasi-continuity with respect to $\capK_{T}$ for bounded psh function which is \emph{stronger} than the usual one for general psh functions with respect to $\capK$ (see \cite{Bedford_Taylor_82}). We refer to Theorem \ref{th-capconvergedecreasing} as a (uniform) strong quasi-continuity of bounded psh functions. 

\proof  We follow ideas presented in \cite[Proposition 1.12]{Kolodziej05}. Let $K  \Subset U$. Let $T_l$ be of bi-dimension $(m,m)$. We will prove that 
\begin{align}\label{limit-capKT}
\int_K (u_k - u) \ddc v_1 \wedge \cdots \wedge \ddc v_m \wedge T_l \to 0
\end{align}
uniformly in psh functions $0 \le v_1, \ldots, v_m \le 1$ and in $l$. The desired assertion concerning the uniform convergence in $\capK_{T_l}$ is a direct consequence of (\ref{limit-capKT}).  

By Hartog's lemma and the boundedness of $u$, we obtain that  $u_k$ is uniformly bounded in $k$ in compact subsets of $U$.  Since the problem is local and  $u_k,u, v_1, \ldots, v_m$ are uniformly bounded on $U$, we can assume that $U\Subset \C^n$,  $u_k,u, v_1, \ldots, v_m$ are defined on an open neighborhood of $\overline U$ and there exist  a smooth psh function $\psi$ defined on an open neighborhood of $\overline U$ and an open neighborhood $W$ of $\partial U$ such that $K \subset U \backslash W$ and  $u_k=u=v_j=\psi$ on $W$ for every $j,k$. Let 
$$T'_l:= \ddc v_2 \wedge \cdots \wedge \ddc v_m \wedge T_l.$$
 Observe $u_k -u$ is of compact support in some  open set $U_1 \Subset U$ containing $K$. Hence, by integration by parts, we get
\begin{align*}
\int_{U_1} (u_k - u) \ddc v_1 \wedge T'_l  &= - \int_{U_1} d(u_k - u) \wedge \dc v_1 \wedge T'_l \\
&\le \bigg(\int_{U_1} d(u_k - u) \wedge \dc (u_k -u) \wedge T'_l\bigg)^{1/2} \bigg(\int_{U_1} d v_1 \wedge \dc v_1 \wedge T'_l \bigg)^{1/2}
\end{align*}
which is $ \lesssim \bigg(\int_{U_1} d(u_k - u) \wedge \dc (u_k -u) \wedge T'_l\bigg)^{1/2}$ by the Chern-Levine-Nirenberg inequality. Denote by $I$ the integral in the last quantity. We have 
 $$I= - \int_{U_1} (u_k - u) \wedge \ddc (u_k -u) \wedge T' \le \int_{U_1} (u_k - u) \wedge \ddc u \wedge T'_l.$$ 
Applying similar arguments to $v_2,\ldots, v_m$ consecutively and the right-hand side of the last inequality, we obtain that 
\begin{align}\label{ine-th-capconvergedecreasingchjuyenhetvj}
 \int_K (u_k - u) \ddc v_1 \wedge \cdots \wedge  \ddc v_m \wedge T_l \le C \bigg( \int_{U_1} (u_k -u) (\ddc u)^m \wedge T_l\bigg)^{2^{-m}},
 \end{align}
where $C$ is independent of $k$ and $l$  (note that the mass of $T_l$ on compact subsets of $U$ is bounded uniformly in $l$). Let 
$$H_{k,l}:= \int_{U_1} (u_k -u) (\ddc u)^m \wedge T_l.$$ 
By (\ref{ine-th-capconvergedecreasingchjuyenhetvj}), in order to obtain  (\ref{limit-capKT}), it suffices to prove that $H_{k,l}$ converges to $0$ as $k \to \infty$ uniformly in $l$. Suppose that this is not the case. This means that there exists a constant $\epsilon>0$, $(k_s)_s\to \infty$ and $(l_s)_s\to \infty$ such that $H_{k_s, l_s} \ge \epsilon$ for every $s$. However, by Condition $(*)$, we get $(u_{k_s} -u) (\ddc u)^{m} \wedge T_{l_s} \to 0$ as $s \to \infty$. This is contradiction.  Hence, (\ref{limit-capKT}) follows.  

We prove the second desired assertion, let $K \Subset U$ and  $(u_k)_k$ a sequence of smooth psh functions defined on an open neighborhood of $K$ decreasing to $u$.  Let $\epsilon>0$ be a constant. Since $u_k \to u$ in $\capK_{T_l}$ as $k \to \infty$ uniformly in $l$, there is a sequence  $(j_k)_k$ converging to $\infty$ for which 
$$\capK_{T_l} \big(K \cap \{ u_{j_k} > u+ 1/k \},  U\big) \le \epsilon 2^{-k}$$
for every $k,l \in \N^*$. Consequently, for $K_\epsilon:= K \backslash \cup_{k=1}^\infty \{ u_{j_k} > u+ 1/k \}$, we have that $\capK_{T_l}(K\backslash K_\epsilon,U) \le \epsilon$ and $u_{j_k}$ is convergent uniformly on $K_\epsilon$. Hence $u$ is continuous on $K_\epsilon$.

Let $(U_s)_s$ be an increasing exhaustive sequence of relatively compact open subsets of $U$ and  $K_s:= \overline U_s \backslash U_{s-1}$ for $s \ge 1$, where $U_0:= \varnothing$. Observe that $K_l$ is  compact,  $U= \cup_{s=1}^\infty K_s$ and 
\begin{align}\label{inclu-KlcapK}
K_s \cap \overline{ \cup_{s' \ge s+2}K_{s'}} =\varnothing
\end{align}
 for  every $s \ge 1$. By the previous paragraph, there exists a compact subset $K'_s$ of $K_s$ such that $\capK_{T_l}(K_s \backslash K'_s,U) \le \epsilon 2^{-s}$ and $u$ is continuous on $K'_s$. Observe that $K':=\cup_{s=1}^\infty K'_s$ is closed in $U$ and $u$ is continuous on $K'$ because of (\ref{inclu-KlcapK}). We also have $U \backslash K' \subset \cup_{s=1}^\infty (K_s \backslash K'_s)$. Hence $\capK_{T_l}(U\backslash K', U)\le \epsilon$ for every $l$. The proof is finished.
\endproof

As one can expect, the above quasi-continuity of bounded psh functions allows us to treat, to certain extent,  these functions as continuous functions with respect to closed positive currents.

\begin{corollary} \label{cor-convergplurifine}  Let $R_k:= \ddc v_{1k} \wedge \cdots \wedge \ddc v_{m k} \wedge T_k$ and $R:=  \ddc v_{1} \wedge \cdots \wedge  \ddc v_{m} \wedge T$, where $v_{j k}, v_j$ are uniformly bounded psh functions on $U$ and $T_k,T$ closed positive currents of bi-degree $(p,p)$.  Let $u$ be a bounded psh function on $U$ and $\chi$ a continuous function on $\R$.  Assume that $R_k \to R$ as $ k\to \infty$ on $U$ and $(T_k)_k$ satisfies Condition $(*)$.
 Then we have 
$$\chi(u) R_k \to \chi(u) R$$
 as $k \to \infty$. In particular, the last convergence holds when $T_k= T$ for every $k$ or $T_k = \ddc w_k\wedge S$, $T= \ddc w\wedge S$, where  $S$ is a closed positive current,  $w$ is a psh function locally integrable with respect to $S$ and $w_k$ is a psh function converging to $w$ in $L^1_{loc}$ as $k \to \infty$ so that $w_k \ge w$ for every $k$.   
\end{corollary}

\proof 
The problem is local. Hence we can assume $U$ is relatively compact in $\C^n$.  Since $u$ is bounded, using Theorem \ref{th-capconvergedecreasing}, we have that  $u$ is uniformly quasi-continuous with respect to the family $\capK_{T_k}$ with $k \in \N$. This means that given $\epsilon>0$, we can find an open subset $U' $ of $U$ such that $\capK_{T_k}(U',U) < \epsilon$ and $u|_{U\backslash U'}$ is continuous.   Let $\tilde{u}$ be a  bounded continuous function on $U$ extending $u|_{U\backslash U'}$ (see \cite[Theorem 20.4]{Rudin}). We have $\chi(\tilde{u}) R_k \to \chi(\tilde{u}) R$ because $\chi, \tilde{u}$ are continuous. 
Moreover, 
$$\big\|\big(\chi(\tilde{u})- \chi(u)\big)R_k\big \| \lesssim \|  R_k\|_{U \backslash U'} \le \capK_{T_k}(U \backslash U', U)  <\epsilon$$
(we used here the boundedness of $U$) and a similar estimate also holds for  $\big(\chi(\tilde{u})- \chi(u)\big)R$. The desired assertion then follows. This finishes the proof.   
\endproof

The following result generalizes well-known convergence properties of Monge-Amp\`ere operators in \cite{Bedford_Taylor_82}.

\begin{theorem} \label{th-increasingsequenceMa} Let $U \subset \C^n$ be an open set.  Let $(T_k)_k$ be a sequence of closed positive currents satisfying Condition $(*)$ so that   $T_k$ converges to  a closed positive current $T$ on $U$ as $k \to \infty$. Let $u_j$ be a locally bounded psh function on $U$ for $1 \le j \le m$. Let $(u_{jk})_{k \in \N}$ be a sequence of locally bounded psh functions converging to $u_j$ in $L^1_{loc}$ as $k \to \infty$. Then, the convergence 
\begin{align}\label{limit-increseinpointwise}
u_{1k} \ddc u_{2k} \wedge \cdots \wedge \ddc u_{mk} \wedge T_k \to u_1 \ddc u_{2} \wedge \cdots \wedge \ddc u_{m} \wedge T
\end{align}
as $k \to \infty$ holds provided that one of the following conditions is fulfilled: 

$(i)$ $u_{jk}(x) \nearrow u_j(x)$  for every $x\in U$ as $k \to \infty$,  

$(ii)$  $u_{jk}(x) \nearrow u_j(x)$ for almost everywhere  $x \in U$ (with respect to the Lebesgue measure) and $T$ has no mass on pluripolar sets,

$(iii)$  $u_{jk} \ge u_j$ for every $j,k$.
\end{theorem}

\proof Given that we already have a uniform strong quasi-continuity for bounded psh functions, the desired result can be deduced without difficulty from proofs of classical results on the convergence of Monge-Amp\`ere operators, for example, see \cite{Kolodziej05}. 

We present here a proof of Theorem \ref{th-increasingsequenceMa} for the readers' convenience. First of all,  observe that if $u_{jk} \nearrow u_j$ almost everywhere then, we have  $u_{jk} \le u_{j (k+1)} \le  u_j$ pointwise on $U$ and the set $\{x \in U:  u_j(x) \not = \lim_{k\to \infty} u_{jk}(x) \}$ is pluripolar. By the localization principle (\cite[Page 7]{Kolodziej05}), we can assume that $u_{jk}, u_j$ are all equal to some smooth psh function $\psi$ outside some set $K \Subset U$ on $U$.   
We prove $(i), (ii)$ simultaneously.   Let 
$$S_{jk}:= \ddc u_{jk} \wedge \cdots \wedge \ddc u_{mk} \wedge T_k, \quad S_{j}:= \ddc u_{j} \wedge \cdots \wedge \ddc u_{m} \wedge T.$$
We prove by induction in $j$ that
\begin{align}\label{ine-increasingsequenceMasosanhRjinfty}
u_{(j-1)k} S_{jk} \to  u_{(j-1)} S_j
\end{align}
 $k$ and  for every $2 \le j \le m+1$ (by convention we put $S_{(m+1)k}:= T_k$ and  $S_{m+1}:=  T$).   The claim is true for $j=m+1$. Suppose that it holds for $(j+1)$. We need to prove it for $j$. Let $R_{j \infty}$ be a limit current of  $u_{(j-1)k} S_{jk}$ as $k \to \infty$.  By induction hypothesis (\ref{ine-increasingsequenceMasosanhRjinfty}) for $(j+1)$ instead of $j$,  $S_{j k} \to S_{j}$ as $k \to \infty$. This combined with the fact that the sequence $(u_{jk})_k$ converges in $L^1_{loc}$ to $u_j$ gives
\begin{align}\label{ine-increasingsequenceMalimsup2}
R_{j \infty} \le u_{j-1} S_{j}
\end{align}
(one can see \cite[Proposition 3.2]{Fornaess_Sibony}).  On the other hand, since $(u_{jk})_k$ is increasing, using Corollary \ref{cor-convergplurifine},  we obtain
\begin{align*} 
\liminf_{k \to \infty} u_{(j-1)k} S_{jk} \ge \liminf_{k \to \infty} u_{(j-1)s} S_{jk}= u_{(j-1)s} S_j
\end{align*}
  for every $s \in \N$. Letting $s \to \infty$ in the last inequality gives 
\begin{align}\label{ine--increasingsequenceMaiminf2}
R_{j \infty}\ge (\lim_{s \to\infty}u_{(j-1)s}) S_j= u_{j-1}S_j + (\lim_{s \to\infty}u_{(j-1)s}- u_{j-1}) S_j.
\end{align}
Recall that the set of  $x \in U$ with $u_{j-1}(x)> \lim_{s \to\infty}u_{(j-1)s}(x) $ is empty in the setting of $(i)$ and is  a pluripolar set in the setting of $(ii)$.   Hence (\ref{ine-increasingsequenceMasosanhRjinfty}) follows from Lemma \ref{le-capTK}, (\ref{ine--increasingsequenceMaiminf2}) and (\ref{ine-increasingsequenceMalimsup2}). We have proved $(i)$ and $(ii)$.  We prove $(iii)$ by similar induction.  The proof is finished.  
\endproof

\begin{remark}\label{re-th-increasingsequenceMa} By the above proof and  Lemma \ref{le-capTK}, Property $(ii)$ of Theorem \ref{th-increasingsequenceMa} still holds if instead of requiring $T$ has no mass on pluripolar sets, we assume the following two conditions:

$(i)$  $T$ has no mass on $A_j:= \{x \in U:  u_j(x) \not = \lim_{k\to \infty} u_{jk}(x) \}$ for every $1 \le j \le m$ and,

$(ii)$  the set $A_j$ is locally complete pluripolar for every $j$. 
\end{remark}

Just by replacing the usual quasi-continuity of psh functions by the stronger one given in  Theorem \ref{th-capconvergedecreasing} for bounded psh functions, we immediately obtain  results similar to those in \cite{BT_fine_87}. We state here  results we will use later. 

\begin{lemma}\label{le-hoitutrentapstreongplurfineopen} (similar to \cite[Lemma 4.1]{BT_fine_87})  Let $U$ be an open subset in $\C^n$.   Let $T$ be a closed positive current on $U$ and $u_j, u_{jk}, u'_j, u'_{jk}$ bounded psh functions on $U$  for $k \in \N$  and $1 \le j \le m$, where $m \in \N$. Let $q \in \N^*$ and  $v_j, v'_j$ bounded psh functions on $U$ for $1 \le j \le q$. Put $W:= \cap_{j=1}^q \{v_j > v'_j\}$.  Assume that  
$$R_k:= \ddc u_{1k} \wedge \cdots \wedge \ddc u_{mk} \wedge T \to R:=\ddc u_{1} \wedge \cdots \wedge \ddc u_{m} \wedge T$$
and 
$$R'_k:= \ddc u'_{1k} \wedge \cdots \wedge \ddc u'_{mk} \wedge T \to R':=\ddc u'_{1} \wedge \cdots \wedge \ddc u'_{m} \wedge T$$
as $k \to \infty$ and 
\begin{align}\label{eq-!WRkRkphaySec1} 
\bold{1}_{W} R_k=\bold{1}_{W} R'_k
 \end{align}
 for every $k$.  Then we have $\bold{1}_{W} R=\bold{1}_{W} R'.$
\end{lemma}

\proof  The problem is clear if $W$ is open, for example, when $v_j$ is continuous for $1 \le j \le q$. In the general  case, we will use the strong quasi-continuity to modify $v_j$. Since the problem is local, we can assume that $U$ is bounded.   Let $\epsilon>0$ be a constant.   By Theorem \ref{th-capconvergedecreasing}, we can find bounded continuous functions $\tilde{v}_j$ on $U$ such that   $\capK_T(\{\tilde{v}_j\not = v_j\},U)< \epsilon$.  Put $\tilde{W}:= \cap_{j=1}^q \{\tilde{v}_j > v'_j\}$ which is an open set.  The choice of $\tilde{v}_j$  combined with the definition of $\capK_T$ yields that 
$$\| \bold{1}_{W} R-\bold{1}_{\tilde{W}} R \|_U \le \epsilon, \quad \| \bold{1}_{W} R_k-\bold{1}_{\tilde{W}} R_k \|_U \le \epsilon.$$ 
We also have similar estimates for $R', R'_k$.  By this and (\ref{eq-!WRkRkphaySec1}),  we get  $ \| \bold{1}_{\tilde{W}} R_k-\bold{1}_{\tilde{W}} R'_k \|_U \le 2 \epsilon$. This combined with the fact that $\tilde{W}$ is open yields that 
$\| \bold{1}_{\tilde{W}} R-\bold{1}_{\tilde{W}} R' \|_U \le 2 \epsilon$. Thus $\| \bold{1}_{W} R-\bold{1}_{W} R' \|_U \le 4 \epsilon$ for every $\epsilon$. The desired equality follows. This finishes the proof.
\endproof

\begin{theorem} \label{th-equalityMAonstrongerplurifinetopo} Let $U$ be an open subset in $\C^n$.   Let $T$ be a closed positive current on $U$ and $u_{j}, u'_j$ bounded psh functions on $U$ for $1 \le j \le m$, where $m \in \N$. Let $v_j,v'_j$ be bounded psh functions on $U$ for $1 \le j \le q$. Assume that $u_{j}= u'_{j}$ on $W:=\cap_{j=1}^q\{v_j> v'_j\}$ for $1 \le j \le m$. Then we have 
\begin{align}\label{eq-MAequaonstrongplurifine}
\bold{1}_{W} \ddc u_{1} \wedge \cdots \wedge \ddc u_{m} \wedge T=\bold{1}_{W} \ddc u'_{1} \wedge \cdots \wedge  \ddc u'_{m} \wedge T.
\end{align}  
\end{theorem}

\proof We give here a complete proof for the readers' convenience.  Let $\epsilon>0$ be a constant. Put $u''_j:= \max\{u_j, u'_j-\epsilon\}$ and $\tilde{W}:= \cap_{j=1}^m \{u_j > u'_j- \epsilon\}$. By hypothesis, $W \subset \tilde{W}$.  
We will prove that 
\begin{align}\label{eq-MAequaonstrongplurifinestep1new}
\bold{1}_{ \tilde{W}}  \ddc u_{1} \wedge \cdots \wedge \ddc u_{m} \wedge T=\bold{1}_{ \tilde{W}}  \ddc u''_{1} \wedge \cdots \wedge \ddc u''_{m} \wedge T.
\end{align} 
Since the problem is local, we can assume there is a sequence of uniformly bounded smooth psh functions $(u_{j k})_k$ decreasing to $u_j$ for $1 \le j \le m$.  Since $\tilde{W}_{k}:= \{u_{jk} >u'_j- \epsilon\}$ is open, we have 
$$\bold{1}_{ \tilde{W}_k} \ddc u_{1k} \wedge \cdots \wedge \ddc u_{mk}\wedge T=\bold{1}_{ \tilde{W}_k} \ddc \max\{u_{1k}, u'_j -\epsilon\} \wedge \cdots \wedge \ddc \{u_{mk}, u'_j -\epsilon\}\wedge T .$$
This together with the inclusion $\tilde{W} \subset \tilde{W}_k$ gives
$$\bold{1}_{ \tilde{W}} \ddc u_{1k} \wedge \cdots \wedge \ddc u_{mk}\wedge T=\bold{1}_{ \tilde{W}} \ddc \max\{u_{1k}, u'_j -\epsilon\} \wedge \cdots \wedge \ddc \{u_{mk}, u'_j -\epsilon\}\wedge T .$$
Using this and Lemma \ref{le-hoitutrentapstreongplurfineopen}, we obtain (\ref{eq-MAequaonstrongplurifinestep1new}) by considering $k \to \infty$. In particular, we get
$$\bold{1}_{ W}  \ddc u_{1} \wedge \cdots \wedge \ddc u_{m} \wedge T=\bold{1}_{W}  \ddc u''_{1} \wedge \cdots \wedge \ddc u''_{m} \wedge T.$$
Letting $\epsilon \to 0$ and using Lemma \ref{le-hoitutrentapstreongplurfineopen} again gives 
$$\bold{1}_{ W}  \ddc u_{1} \wedge \cdots \wedge \ddc u_{m} \wedge T=\bold{1}_{W}  \ddc \max\{u_1, u'_1\} \wedge \cdots \wedge \ddc \max\{u_m, u'_m\} \wedge T.$$
The last equality still holds if we replace $u_j$ in the left-hand side by $u'_j$ by using similar arguments. So the desired equality follows.  The proof is finished. 
\endproof

\begin{remark} \label{re-quasipsh} Recall that a quasi-psh function $u$ on $U$ is, by definition, locally the sum of a psh function and a smooth one. We can check that  results presented above have their analogues for quasi-psh functions.
\end{remark}

\section{Relative non-pluripolar product} \label{sec-def-gene-nonpluri}

Let $X$ be a complex manifold of dimension $n$ and  $T, T_1, \ldots, T_m$  closed positive currents on $X$ such that $T_j$ is of bi-degree $(1,1)$ for $1 \le j \le m$.  
Let $U$ be a local chart of $X$ such that $T_j= \ddc u_j$ on $U$ for $1 \le j \le m$, where   $u_1, \ldots, u_m$ are psh functions on $U$. Let $k \in \N$ and $u_{jk}:= \max\{u_j, -k\}$ which is a locally bounded psh function.  Put $R_k:= \ddc u_{1k} \wedge \cdots \wedge \ddc u_{mk} \wedge T$. By Theorem \ref{th-equalityMAonstrongerplurifinetopo} and the fact that $\{u_j > -k\}= \{u_{jk} > -k\}$, we have 
\begin{align}\label{eq-Rkbangnhautrentapulonhonkgeneralized}
\bold{1}_{\cap_{j=1}^m \{u_j > -k\}} R_k= \bold{1}_{\cap_{j=1}^m \{u_j > -k\}} R_{l}
\end{align}
for every $l \ge k$. As in the case of the usual non-pluripolar products, we have the following basic observation.

\begin{lemma} \label{le-basicprononpluri} Assume that we have
 \begin{align}\label{ine-dieukiendecogeneralizedmass}
\sup_{k\in \N} \| \bold{1}_{\cap_{j=1}^m \{u_j > -k\}} R_k \|_K < \infty
\end{align}
for every compact $K$ of $U$. Then the limit current
\begin{align}\label{eq-defcuageneralnonpluripolar}
R:= \lim_{k\to \infty }\bold{1}_{\cap_{j=1}^m \{u_j > -k\}} R_k
\end{align}
is well-defined and  for every Borel  form $\Phi$ with bounded coefficients on $U$ such that $\supp \Phi \Subset U$, we have 
\begin{align}\label{eq-le-basicprononpluri}
\langle R, \Phi\rangle = \lim_{k \to \infty} \langle \bold{1}_{\cap_{j=1}^m \{u_j > -k\}} R_k, \Phi\rangle.
\end{align}
Consequently, there holds 
$$\bold{1}_{\cap_{j=1}^m \{u_j > -k\}} R= \bold{1}_{\cap_{j=1}^m \{u_j > -k\}} R_k, \quad \bold{1}_{\cup_{j=1}^m \{u_j = -\infty\}} R= 0.$$ 
\end{lemma}

\proof  By (\ref{eq-Rkbangnhautrentapulonhonkgeneralized}), we have 
\begin{align*}
\bold{1}_{\cap_{j=1}^m \{u_j > -l\}} R_l &= \bold{1}_{\cap_{j=1}^m \{ u_j > 0\}} R_l+ \sum_{k=1}^l \bold{1}_{\cap_{j=1}^m \{-k+1 \ge  u_j > -k\}} R_l \\
&=\bold{1}_{\cap_{j=1}^m \{u_j > 0\}} R_0+ \sum_{k=1}^l \bold{1}_{\cap_{j=1}^m \{-k+1 \ge u_j > -k \}} R_k.
\end{align*}
This combined with (\ref{ine-dieukiendecogeneralizedmass}) tells us that  the mass on a fixed compact of $U$ of the current
$$\bold{1}_{\cap_{j=1}^m \{-k \ge  u_j > -l  \}} R_l= \sum_{k'=k}^l \bold{1}_{\cap_{j=1}^m \{-k'+1 \ge u_j > -k' \}} R_k$$
converging to $0$ as $l\ge k\to \infty$.  We deduce that $\lim_{k\to \infty }\bold{1}_{\cap_{j=1}^m \{u_j > -k\}} R_k$ exists and is denoted by $R$. 

Since $R, R_k$ are positive, for every continuous form $\Phi$ of compact support in $U$, we have $\langle R, \Phi\rangle = \lim_{k \to \infty} \langle R_k, \Phi\rangle$.  Let $\Phi$ be a Borel form on $U$ such that its coefficients are bounded on $U$ and $\supp \Phi \Subset U$.  Let $K$ be a compact of $U$ containing $\supp \Phi$ and $U_1 \supset K$ a relatively compact open subset of $U$. Let $\epsilon>0$ be a constant. Let $k_0$ be a positive integer such that 
\begin{align}\label{ine-uocliongkolusintheoremPhi[hi'0}
\|\bold{1}_{\cap_{j=1}^m \{-k \ge  u_j > -l \}} R_l\|_{U_1} \le \epsilon
\end{align}
for every $l\ge k \ge k_0$. By Lusin's theorem, there exists a continuous form $\Phi'$ compactly supported on $U_1$ such that  
\begin{align}\label{ine-uocliongkolusintheoremPhi[hi'}
\|\bold{1}_{\cap_{j=1}^m \{u_j > -k_0\}} R_{k_0} \|_{\{x \in U_1: \Phi'(x) \not= \Phi(x)\}}\le \epsilon, \quad \|R\|_{\{x \in U_1: \Phi'(x) \not= \Phi(x)\}}\le \epsilon.
\end{align}
Using (\ref{ine-uocliongkolusintheoremPhi[hi'0}) and (\ref{ine-uocliongkolusintheoremPhi[hi'}) gives
$$|\langle |\bold{1}_{\cap_{j=1}^m \{u_j > -k\}} R_k, \Phi\rangle-\langle R_k, \Phi'\rangle| \lesssim 2 \epsilon, \quad |\langle R, \Phi\rangle-\langle R, \Phi'\rangle| \lesssim 2 \epsilon.$$
This combined with the fact that $|\langle |\bold{1}_{\cap_{j=1}^m \{u_j > -k\}} R_k, \Phi'\rangle \to  \langle R, \Phi'\rangle$ gives 
$$ \big| \lim_{k \to \infty}\langle  |\bold{1}_{\cap_{j=1}^m \{u_j > -k\}} R_k, \Phi\rangle -\langle R, \Phi\rangle \big| \lesssim 2 \epsilon.$$
Letting $\epsilon \to 0$ gives (\ref{eq-le-basicprononpluri}). For the other equalities, one just needs to apply (\ref{eq-le-basicprononpluri}) to suitable $\Phi$. This finishes the proof.
\endproof

\begin{lemma} \label{le-globalpotentialnonpluri} Let $T_j= \ddc \tilde{u}_j+ \theta_j$ for $1 \le j \le m$, where $\theta_j$ is a smooth $(1,1)$-form and $\tilde{u}_j$ is  $\theta_j$-psh on $X$. Let 
$$\tilde{u}_{jk}:= \max\{\tilde{u}_j, -k\}, \quad \tilde{R}_k:= \wedge_{j=1}^m (\ddc \tilde{u}_{jk}+ \theta_j) \wedge T.$$
Then the following two properties hold:

$(i)$ (\ref{ine-dieukiendecogeneralizedmass}) holds for every small enough local chart $U$ if and only if we have that  for every compact $K$ of $X$, 
\begin{align}\label{ine-dieukiendecogeneralizedmassglobal}
\sup_{k \in \N} \| \bold{1}_{\cap_{j=1}^m \{\tilde{u}_j >-k\}} \tilde{R}_k\|_K < \infty.
\end{align}
In this case, if $\tilde{R}:= \lim_{k \to \infty}  \bold{1}_{\cap_{j=1}^m \{\tilde{u}_j >-k\}}\tilde{R}_k$, then  $\tilde{R}=R$ on $U$,

$(ii)$  $\bold{1}_{\cap_{j=1}^m \{\tilde{u}_j >-k\}} \tilde{R}_k$ is a positive current.
\end{lemma}

\proof Firstly observe that $\tilde{R}_k$ is a current of order $0$ and of bounded mass on compact subsets of $X$.  Let 
$$\tilde{B}_k:=\cap_{j=1}^m \{\tilde{u}_j >-k\}.$$
Assume now (\ref{ine-dieukiendecogeneralizedmassglobal}). This means   $\| \bold{1}_{\tilde{B}_k} \tilde{R}_k\|_K$ is uniformly bounded for every compact $K$.     By Remark \ref{re-quasipsh}, we have $\bold{1}_{\tilde{B}_k} \tilde{R}_k= \bold{1}_{\tilde{B}_k}\tilde{R}_l$ for every $l\ge k$. Decompose 
$$\tilde{B}_l= \tilde{B}_0 \cup_{k=1}^l \cap_{j=1}^m \{-k+1 \ge \tilde{u}_j >-k\}$$
which is a disjoint union. Hence, we get
\begin{align*}
\|\bold{1}_{\tilde{B}_l} \tilde{R}_l\| = \|\bold{1}_{\tilde{B}_0} \tilde{R}_l\|+ \sum_{k=1}^l \|\bold{1}_{\cap_{j=1}^m \{-k+1 \ge  \tilde{u}_j > -k\}} \tilde{R}_l\| =\|\bold{1}_{\tilde{B}_0} \tilde{R}_0\|+ \sum_{k=1}^l \|\bold{1}_{\cap_{j=1}^m \{-k+1 \ge \tilde{u}_j > -k \}} \tilde{R}_k\|,
\end{align*}
where the masses are measured on some compact $K$ in $X$. We deduce that the condition $\| \bold{1}_{\tilde{B}_k} \tilde{R}_k\|_K$ is uniformly bounded is equivalent to that 
\begin{align}\label{conver-dktuonduongRnganonpluri}
 \|\bold{1}_{\cap_{j=1}^m \{-k \ge \tilde{u}_j > -l\}} \tilde{R}_l\|_K = \sum_{k'=k}^l \|\bold{1}_{\cap_{j=1}^m \{-k'+1 \ge \tilde{u}_j > -k' \}} \tilde{R}_k\|_K \to 0
\end{align}
 as $l \ge k \to \infty$. Hence $\bold{1}_{\tilde{B}_k} \tilde{R}_k$ converges to a current $\tilde{R}$ and moreover we have 
$$\langle \tilde{R}, \Phi\rangle = \lim_{k \to \infty} \langle \bold{1}_{\tilde{B}_k} \tilde{R}_k, \Phi\rangle$$
for every Borel bounded form $\Phi$ of compact support in $X$ as in the proof of Lemma \ref{le-basicprononpluri}. Consequently, we get
\begin{align} \label{eq-RngabangRngak}
\bold{1}_{\tilde{B}_k} \tilde{R}= \bold{1}_{\tilde{B}_k} \tilde{R}_k, \quad \bold{1}_{\cup_{j=1}^m \{\tilde{u}_j = -\infty\}} \tilde{R}= 0.
\end{align}

Let $U, u_j, u_{jk}, R, R_k$ be as above. We will show that (\ref{ine-dieukiendecogeneralizedmass}) is satisfied and  $\tilde{R}= R$ on $U$.  Let 
$$B_k:=\cap_{j=1}^m \{u_j >-k\}.$$
Observe that $u_j=\tilde{u}_j+ \tau_j$ for some smooth function $\tau_j$ on $U$ with $\ddc \tau_j = \theta_j$. By shrinking $U$, we can assume that $\tau_j$ is bounded on $U$ and let $c_0$ be an integer greater than  $\sum_{j=1}^m \|\tau_j\|_{L^\infty}$. 
We have 
$$\tilde{u}_{j k}+ \tau_j= \max\{ \tilde{u}_j+ \tau_j, -k + \tau_j\}$$
 which is equal to $\max\{ \tilde{u}_j+ \tau_j, -k \}= u_{jk}$ on the set  $\{\tilde{u}_j > -k + c_0\}$. It follows that 
\begin{align}\label{eq-bieudienRjkquaRj}
 \bold{1}_{\tilde{B}_{k-c_0}} \tilde{R}_k =  \bold{1}_{\tilde{B}_{k-c_0}} \big(\wedge_{j =1}^m \ddc u_{jk} \wedge T) = \bold{1}_{\tilde{B}_{k-c_0}}  R_k.
 \end{align}
 This together with the inclusions  $B_{k-2 c_0} \subset \tilde{B}_{k-c_0} \subset B_k$ give
\begin{align}\label{eq-bieudienRjkquaRj23}
 \bold{1}_{B_{k-2c_0}}  R_{k-2 c_0} \le  \bold{1}_{\tilde{B}_{k-c_0}} \tilde{R}_k= \bold{1}_{\tilde{B}_{k-c_0}} \tilde{R}_{k-c_0} \le  \bold{1}_{B_{k}} R_k.
 \end{align}
Hence (\ref{ine-dieukiendecogeneralizedmass}) follows. We also deduce from this and (\ref{eq-RngabangRngak}) that $\tilde{R}=R$. Conversely, if  (\ref{ine-dieukiendecogeneralizedmass}) holds for every $U$, then,  by (\ref{eq-bieudienRjkquaRj23}), the claim (\ref{ine-dieukiendecogeneralizedmassglobal}) holds. Hence $(i)$ follows.  By (\ref{eq-bieudienRjkquaRj}), we have 
$$\bold{1}_{\tilde{B}_{k}} \tilde{R}_k= \bold{1}_{\tilde{B}_{k}} \tilde{R}_{k+ c_0}= \bold{1}_{\tilde{B}_{k}} R_{k+c_0} \ge 0.$$
Thus $(ii)$ follows.  The proof is finished.
\endproof

Lemma \ref{le-globalpotentialnonpluri} applied to $U$ in place of $X$ implies that the condition (\ref{ine-dieukiendecogeneralizedmass}) is independent of the choice of $u_j$ and so is the limit $R$ above.  As a result, if  (\ref{ine-dieukiendecogeneralizedmass}) holds for every small enough local chart $U$ as above, 
then we obtain a positive current $R$ globally defined on $X$ given locally by (\ref{eq-defcuageneralnonpluripolar}). This current is equal to $\tilde{R}$ by Lemma \ref{le-globalpotentialnonpluri} again.   

\begin{definition} We say that  the non-pluripolar product relative to $T$ of $T_1, \ldots, T_m$ is \emph{well-defined} if (\ref{ine-dieukiendecogeneralizedmass}) holds for every small enough local chart $U$ of $X$, or equivalently, (\ref{ine-dieukiendecogeneralizedmassglobal}) holds. In this case, \emph{the non-pluripolar product relative to $T$}  of $T_1, \ldots, T_m$, which is denoted by $\langle T_{1} \wedge \cdots \wedge T_{m} \dot{\wedge} T \rangle$,  is  defined to be  the current $\tilde{R}$ in Lemma \ref{le-globalpotentialnonpluri}.
\end{definition}

When $T$ is the current of integration along $X$,  we simply write $\langle T_1 \wedge \cdots \wedge T_m \rangle$ for the non-pluripolar product relative to $T$ of $T_1, \ldots, T_m$. One can see that this is exactly the usual non-pluripolar product of $T_1, \ldots, T_m$ defined in \cite{BT_fine_87,GZ-weighted,BEGZ}.  We note that  in general the current  $\langle T_1 \wedge \cdots \wedge T_m \dot{\wedge} T \rangle$ can have mass on pluripolar sets in $\cap_{j=1}^m \{u_j >-\infty\}$, see,  however, Property $(iii)$ of Proposition \ref{pro-sublinearnonpluripolar} below.  Arguing as in the proof of  \cite[Proposition 1.6]{BEGZ}, we obtain the following result.

\begin{lemma} \label{le-compactKahlernonpluripolar} Assume that  $X$ is a compact K\"ahler manifold. Then,  $\langle T_1 \wedge \cdots \wedge T_m \dot{\wedge} T \rangle$ is well-defined.
\end{lemma}

For a closed positive $(1,1)$-current $R$,  recall that the \emph{polar locus} $I_R$ of $R$ is  the set of points where potentials of $R$ are equal to $-\infty$.  Note that by Lemma \ref{le-basicprononpluri}, the current  $\langle  T_1 \wedge \cdots \wedge T_m \dot{\wedge} T\rangle$ has no mass on $\cup_{j=1}^m I_{T_j}$.  We collect here some more basic properties of relative non-pluripolar products. 

\begin{proposition}\label{pro-sublinearnonpluripolar} 
$(i)$ The product $ \langle  T_1 \wedge \cdots \wedge T_m \dot{\wedge} T\rangle $ is symmetric with respect to $T_1, \ldots, T_m$. 

$(ii)$   Given a positive real number $\lambda$, we have $\langle (\lambda T_1) \wedge T_2 \wedge \cdots \wedge T_m \dot{\wedge} T \rangle = \lambda \langle T_1 \wedge T_2 \wedge \cdots \wedge T_m \dot{\wedge} T\rangle $.

$(iii)$ Given a complete pluripolar set $A$ such that $T$ has no mass on $A$, then $\langle  T_1 \wedge T_2 \wedge \cdots \wedge T_m \dot{\wedge} T\rangle $ also has no mass on $A$.

$(iv)$  Let $T'_1$ be  a closed positive $(1,1)$-current on $X$ and $T_j,T$ as above. Assume that $\langle \wedge_{j=1}^m T_j \dot{\wedge} T\rangle$, $\langle T'_1 \wedge \wedge_{j=2}^m T_j \dot{\wedge} T\rangle$ are well-defined. Then $\langle (T_1+T'_1) \wedge \wedge_{j=2}^m T_j \dot{\wedge} T\rangle$ is also well-defined and satisfies
\begin{align}\label{ine-convexnonpluripoarTjTphayj}
\big\langle  (T_1+T'_1) \wedge \wedge_{j=2}^m T_j \dot{\wedge} T\big\rangle \le \langle T_1 \wedge \wedge_{j=2}^m T_j \dot{\wedge}T\rangle+ \langle T'_1 \wedge  \wedge_{j=2}^m  T_j \dot{\wedge} T\rangle.
\end{align}
The equality occurs if  $T$ has no mass on $I_{T_1} \cup I_{T'_1}$. 

$(v)$ Let $1 \le l \le m$ be an integer.  Let $T''_j$ be  a closed positive $(1,1)$-current on $X$ and $T_j,T$ as above for $1 \le j \le l$. Assume that $T''_j \ge T_j$ for every $1 \le j \le l$ and $T$ has no mass on $\cup_{j=1}^l I_{T''_j-T_j}$. Then, we have 
$$\langle \wedge_{j=1}^l T''_j \wedge \wedge_{j=l+1}^m T_j \dot{\wedge} T \rangle \ge \langle \wedge_{j=1}^m T_j \dot{\wedge} T \rangle$$
provided that the left-hand side is well-defined. 

$(vi)$ Let $1 \le l \le m$ be an integer. Assume $R:= \langle \wedge_{j=l+1}^m T_j \dot{\wedge} T \rangle$ and $\langle \wedge_{j=1}^l T_j \dot{\wedge} R \rangle$  are well-defined. Then, we have $\langle \wedge_{j=1}^m T_j \dot{\wedge} T \rangle = \langle \wedge_{j=1}^l T_j \dot{\wedge} R \rangle$.

$(vii)$ Let $A$ be a complete pluripolar set. Assume that $\langle  T_1 \wedge T_2 \wedge \cdots \wedge T_m \dot{\wedge} T\rangle $  is well-defined. Then we have 
$$\bold{1}_{X \backslash A}\langle  T_1 \wedge T_2 \wedge \cdots \wedge T_m \dot{\wedge} T\rangle= \big\langle  T_1 \wedge T_2 \wedge \cdots \wedge T_m \dot{\wedge}(\bold{1}_{X \backslash A} T)\big\rangle.$$
In particular,  the equality
$$\langle \wedge_{j=1}^m T_j \dot{\wedge} T \rangle = \langle \wedge_{j=1}^m T_j \dot{\wedge} T' \rangle$$
holds, where $T':= \bold{1}_{X \backslash \cup_{j=1}^m I_{T_j}} T$.

\end{proposition}


\proof Properties $(i), (ii)$ are clear from the definition and the proof of Lemma \ref{le-basicprononpluri}. We now check $(iii)$.  Let $R:= \langle \wedge_{j=1}^m T_j \dot{\wedge} T\rangle$. Since $T$ has no mass on the complete pluripolar set $A$, using Lemma \ref{le-capTK} gives $\bold{1}_{A}\wedge_{j=1}^m \ddc u_{jk} \wedge T=0$ for every $k$. Using this and Lemma \ref{le-basicprononpluri},  we deduce that
$$\bold{1}_{A} \bold{1}_{ \cap_{j=1}^m \{u_j> -k\}} R= \bold{1}_{A} \bold{1}_{ \cap_{j=1}^m \{u_j> -k\}} \wedge_{j=1}^m \ddc u_{jk} \wedge T =0.$$
Hence $R=0$ on $A$ and $(iii)$ follows. 

We prove $(iv)$. We work  on a small local chart $U$. Write $T'_j= \ddc u'_j$, $u'_{jk}:= \max\{u'_j, -k\}$. Recall $T_j= \ddc u_j$. We can assume $u_j, u'_j \le 0$. Let $R':= \langle T'_1 \wedge \wedge_{j=2}^m T_j \dot{\wedge} T\rangle$ and $R'':= \langle(T_1+T'_1) \wedge \wedge_{j=2}^m T_j \dot{\wedge} T\rangle$. We have 
$$\max\{u_1+ u'_1, -k\}= u_{1k}+ u'_{1k}$$
 on $\{u_1 + u'_1 >-k\}$ because the last set is contained in $\{u_1 > -k\} \cap \{u'_1 > -k \}$. This combined with Lemma \ref{le-basicprononpluri}  yields
$$\bold{1}_{B''_k} R''=\bold{1}_{B''_k} R + \bold{1}_{B''_k} R',$$
where $B''_k: =\{u_1+ u'_1 > -k\} \cap \cap_{j=2}^m \{u_j > -k\}$. Letting $k \to \infty$ in the last equality gives (\ref{ine-convexnonpluripoarTjTphayj}).  
Observe that 
$$\bold{1}_{U \backslash B''_k } R \to \bold{1}_{I_{T_1}\cup I_{T'_1} \cup \, \cup_{j=2}^m I_{T_j}} R$$
 as $k \to \infty$. The last limit is equal to $\bold{1}_{I_{T'_1}} R$ because $R$ has no mass on $\cup_{j=1}^m I_{T_j}$. Moreover since $I_{T_1'}$ is complete pluripolar and $T$ has no mass on $I_{T'_1}$, by $(iii)$, we have that $\bold{1}_{I_{T'_1}} R =0$. Consequently,   $\bold{1}_{U \backslash B''_k } R \to 0$ as $k \to \infty$ and a similar property of $R'$ also holds. Thus we obtain the equality in  (\ref{ine-convexnonpluripoarTjTphayj}) if  $T$ has no mass on $I_{T_1} \cup I_{T'_1}$. 
To get $(v)$, we just need to decompose $T''_j= T_j+ T'_j$ for some closed positive current $T'_j$ and use similar arguments as in the proof of $(iv)$. 
 
We prove $(vi)$. We can assume $u_j \le 0$ for every $1 \le j \le m$.  Let 
$$\psi_{k}:= k^{-1} \max\big\{ \sum_{j=1}^m u_j, -k  \big\}+ 1.$$ 
Observe that $0 \le  \psi_k \le 1$ and 
$$\psi_k=0 \quad \text{on } \, \cup_{j=1}^m \{u_j \le -k\}.$$ 
Note that $R$ has no mass on $I_{T_j}$ for $l+1 \le j \le m$. This combined with $(iii)$ yields that $\langle \wedge_{j=1}^l T_j \dot{\wedge} R \rangle$ gives no mass on $I_{T_j}$ for $1 \le j \le m$.
 Using this and the fact that $\psi_k \nearrow \bold{1}_{X \backslash \cup_{j=1}^m I_{T_j}}$ yields 
\begin{align}  \label{eq-tinhwedgedenTtachtul}
 \langle \wedge_{j=1}^l T_j \dot{\wedge} R \rangle = \lim_{k \to \infty} \psi_k  \langle \wedge_{j=1}^l T_j \dot{\wedge} R \rangle=\lim_{k \to \infty} \psi_k  \wedge_{j=1}^l \ddc u_{jk} \wedge R.
\end{align} 
Now since $\psi_k R = \psi_k \wedge_{j=l+1}^m \ddc u_{jk} \wedge T$ (Lemma \ref{le-basicprononpluri}), we get 
\begin{align}\label{eq-psikRghepvaotrong}
\psi_k  \wedge_{j=1}^l \ddc u_{jk} \wedge R= \psi_k  \wedge_{j=1}^m \ddc u_{jk} \wedge T.
\end{align}
To see why the last equality is true, notice that it is clear if $u_{jk}$'s are smooth and in general, we can  use sequences of smooth psh functions decreasing to $u_{jk}$ for $1 \le j \le l$ and  the convergences of Monge-Amp\`ere operators to obtain (\ref{eq-psikRghepvaotrong}). Combining (\ref{eq-psikRghepvaotrong}) with (\ref{eq-tinhwedgedenTtachtul}) gives the desired assertion.

It remains to prove $(vii)$.  Let $\psi_k$ be as above.  
Let $A= \{\varphi = -\infty\}$, where $\varphi$ is a negative psh function. Define $\psi'_k:=k^{-1} \max\big\{ \varphi+ \sum_{j=1}^m u_j, -k  \big\}+ 1.$ Since $0 \le \psi'_k \le \psi_k$, we get $\{\psi_k= 0\} \subset \{\psi'_k=0\}$.  It follows that
$$\psi'_k \langle  T_1 \wedge T_2 \wedge \cdots \wedge T_m \dot{\wedge} T\rangle= \psi'_k  \wedge_{j=1}^m \ddc u_{jk}\wedge T =\psi'_k  \wedge_{j=1}^m \ddc u_{jk}\wedge (\bold{1}_{X \backslash A} T)$$
because $T = \bold{1}_{X \backslash A} T$ on $\{\psi'_k \not = 0\}$ and the convergence of Monge-Amp\`ere operators. Letting $k \to \infty$ gives the desired assertion.  The proof is finished. 
\endproof

The following result clarifies the relationship between the relative non-pluripolar product and some other known notions of intersection. 

\begin{proposition} \label{pro-tinhnonpluripolarclassically} $(i)$ Let $U, u_j, u_{jk}, R_k$ be as above. Assume that $u_j$ is locally integrable with respect to $T_{j+1} \wedge \cdots \wedge T_m \wedge T$ for $1\le j \le m$ and for every bounded psh function $v$, the current  $v R_k$ converges to $v \, T_1 \wedge \cdots \wedge T_m \wedge T$ on $U$  as $k \to \infty$. Then, the current $\langle T_1 \wedge \cdots \wedge T_m \dot{\wedge} T \rangle $ is well-defined on $U$ and 
\begin{align}\label{eq-gennonpluriclassi}
\langle T_1 \wedge \cdots \wedge T_m \dot{\wedge} T \rangle = \bold{1}_{X \backslash \cup_{j=1}^m I_{T_j}} T_1 \wedge \cdots \wedge T_m \wedge T.
\end{align}
In particular, if $u_1, \ldots, u_{m-1}$ are  locally bounded and $u_m$ is locally integrable with respect to $T$, then we have 
\begin{align} \label{eq-gennonpluriclassi2} 
\langle T_1 \wedge \cdots \wedge T_m \dot{\wedge} T \rangle= \bold{1}_{X\backslash I_{T_m}} T_1 \wedge \cdots \wedge T_m \wedge T.
\end{align}

$(ii)$ If $T$ is of bi-degree $(1,1)$, then we have 
\begin{align}\label{eq-gennonpluriclassi3}
\langle T_1 \wedge \cdots \wedge T_m  \wedge T \rangle= \bold{1}_{X \backslash I_{T}}\langle T_1 \wedge \cdots \wedge T_m \dot{\wedge} T \rangle,
\end{align}
where the left-hand side is the usual non-pluripolar product of $T_1, \ldots, T_m,T$. 
\end{proposition} 

We note that in the above $(i)$, the current $T_j \wedge \cdots \wedge T_m \wedge T$ $(1 \le j \le m)$ is defined inductively as in the classical case (see \cite{Bedford_Taylor_82}).  The assumption of $(i)$ of  Proposition \ref{pro-tinhnonpluripolarclassically} is satisfied in well-known classical contexts, we refer to \cite{Fornaess_Sibony} and references therein for details. We notice also that by $(vii)$ of Proposition \ref{pro-sublinearnonpluripolar}, the right-hand side of (\ref{eq-gennonpluriclassi3}) is equal to $\big\langle T_1 \wedge \cdots \wedge T_m \dot{\wedge}  (\bold{1}_{X \backslash I_{T}}T) \big\rangle$.

\proof We prove $(i)$. The question is local. Hence, we can assume $u_j \le 0$ for every $j$.  Let $\psi_k:= k^{-1} \max\{\sum_{j=1}^m u_j, -k \}+1$. By hypothesis, we get 
\begin{align}\label{converge-paikRr}
\psi_k R_r \to \psi_k T_1 \wedge \cdots \wedge T_m \wedge T
\end{align}
 as $r \to \infty$. Since $\psi_k=0$ on $\{u_j \le -k\}$ and $u_{jr}=u_{jk}$ on $\{u_j > -k \}$  for $r \ge k$, using Theorem \ref{th-equalityMAonstrongerplurifinetopo}, we infer that  $\psi_k R_r =\psi_k R_k.$
Consequently, the mass $\psi_k R_k$ on compact subsets of $U$ is bounded uniformly in $k$. Hence, $\langle T_1 \wedge \cdots \wedge T_m \dot{\wedge} T \rangle $ is well-defined on $U$. The above arguments also show that
$$\psi_k R_r =\psi_k R_k= \psi_k \bold{1}_{\cap_{j=1}^m\{ u_j > -k\}}R_r= \psi_k \bold{1}_{\cap_{j=1}^m\{u_j > -k\}}\langle T_1 \wedge \cdots \wedge T_m \dot{\wedge} T \rangle$$
which is equal to $\psi_k \langle T_1 \wedge \cdots \wedge T_m \dot{\wedge} T \rangle.$ Letting $r \to \infty$ in the last equality and using (\ref{converge-paikRr}) give 
$$\psi_k \, T_1 \wedge \cdots \wedge T_m \wedge T= \psi_k \langle T_1 \wedge \cdots \wedge T_m \dot{\wedge} T \rangle.$$
Letting $k \to \infty$ gives (\ref{eq-gennonpluriclassi}).  To obtain (\ref{eq-gennonpluriclassi2}), we just need to combine  (\ref{eq-gennonpluriclassi}) with Theorem \ref{th-decreasing-unbounded}.

It remain to check $(ii)$. We work locally. Let $u$ be a local potential of $T$ and $u \le 0$. Let $u_k:= \max\{u,-k\}$ and  $\psi'_k:= k^{-1} \max\{u+\sum_{j=1}^m u_j, -k \}+1$. Note that  
$$\psi'_k \langle T_1 \wedge \cdots \wedge T_m  \wedge T \rangle =\psi'_k \wedge_{j=1}^m \ddc u_{jk} \wedge \ddc u_k.$$
Since  $\{\psi_k =0\} \subset \{\psi'_k =0\}$ ($\psi_k \ge \psi'_k\ge 0$), we have
\begin{align*}
\psi'_k \langle T_1 \wedge \cdots \wedge T_m  \dot{\wedge} T \rangle  &=\psi'_k R_k= \psi'_k \wedge_{j=1}^m \ddc u_{jk} \wedge \ddc u\\
&=\psi'_k  \bold{1}_{\{u > -k\}}\wedge_{j=1}^m \ddc u_{jk} \wedge \ddc u \\
&= \psi'_k  \bold{1}_{\{u > -k\}}\wedge_{j=1}^m \ddc u_{jk} \wedge \ddc u_k
\end{align*}
(one can obtain the last equality as a consequence of (\ref{eq-gennonpluriclassi2}) applied to  the case where $T$ is the current of integration along $X$ or alternatively we can use Theorem \ref{th-decreasing-unbounded} directly).
Letting $k \to \infty$ in the last equality gives (\ref{eq-gennonpluriclassi3}). This finishes the proof. 
\endproof

As in the case of the usual non-pluripolar products, the relative non-pluripolar products, if  well-defined,  are closed positive currents as showed by the following result. 

\begin{theorem} Assume that $\langle T_1 \wedge \cdots \wedge T_m \dot{\wedge} T \rangle$ is well-defined. Then $\langle T_{1} \wedge \cdots \wedge T_{m} \dot{\wedge} T \rangle$ is closed. 
\end{theorem}

\proof The proof is  based on ideas from \cite{Sibony_duke,BEGZ}.  We work  on a local chart $U$ as above. By shrinking $U$ and subtracting from $u_j$ a suitable constant, we can assume that $u_j \le 0$.   Let 
$$\psi_k:= k^{-1}\max\{ \sum_{j=1}^m u_j, -k \} + 1.$$
Observe that $\psi_k= 0$ on  $\cup_{j=1}^m \{u_j \le -k\}$ and   $0 \le \psi_k \le 1$ increases to $\bold{1}_{\cap_{j=1}^m \{u_j > -\infty\}}$. Let $g: \R \to \R$ be a non-negative smooth function such that $g(0)=g'(0)=g'(1)=0$ and $g(1)=1$. Let $R_k, R$ be as above.  Since $g(1)=1$ and $g(0)=0$, we get $g(\psi_k)R \to R$ as $k \to \infty$ (recall $R$ has no mass on $\cup_{j=1}^m \{u_j = -\infty\}$). Thus the desired assertion is equivalent to proving that 
\begin{align}\label{eq-limitdgeneralnonpluriproduct}
d R= \lim_{k \to \infty}d(g(\psi_k) R) =0.
\end{align}  
Since $g(\psi_k)=g(0)=0$ on $\cup_{j=1}^m \{u_j \le -k\}$, we have 
$$d(g(\psi_k) R)= d(g(\psi_k) R_k)= g'(\psi_k) d \psi_k \wedge  R_k$$
(see \cite[Lemma 1.9]{BEGZ} or Corollary \ref{cor-convergplurifine} for the second equality).  Let $U_1 \Subset U_2$ be relatively compact open subsets of $U$. By the Cauchy-Schwarz inequality, we have 
$$\| g'(\psi_k) d \psi_k \wedge  R_k\|_{U_1} \lesssim \| d \psi_k \wedge \dc \psi_k \wedge  R_k \|_{U_1} \|g'(\psi_k)^2  R_k \|_{U_1}.$$
Using $ d \psi_k \wedge \dc \psi_k= \ddc \psi_k^2 - \psi_k \ddc \psi_k$ and the Chern-Levine-Nirenberg inequality,   one get 
$$ \| d \psi_k \wedge \dc \psi_k \wedge  R_k \|_{U_1} \lesssim \| \psi_k R_k\|_{U_2} \le  \|R_k\|_{\overline U_2 \cap \cap_{j=1}^m\{ u_j > -k\}} \lesssim 1$$
by (\ref{ine-dieukiendecogeneralizedmass}). On the other hand, using the equality $g'(\psi_k)= g'(0)=0$ on $\{u_j \le -k\}$, we obtain  
$$g'(\psi_k)^2 R_k =g'(\psi_k)^2  \bold{1}_{\cap_{j=1}^m \{u_j > -k\}}R_k = g'(\psi_k)^2  R $$
converging to $0$ as $k \to \infty$ because $g'(\psi_k) \to g'(1)= 0$ on $\cap_{j=1}^m\{u_j >-\infty \}$. Thus  (\ref{eq-limitdgeneralnonpluriproduct}) follows. This finishes the proof.
\endproof

\begin{remark} \label{re-currentinterV} Let $V$ be a smooth submanifold of $X$ and $T$ the current of integration along $V$.  If $V \subset \cup_{j=1}^m I_{T_j}$, then we know that the non-pluripolar product relative to $T$ of $T_1, \ldots, T_m$ is zero. Consider now the case where $V \not \subset \cup_{j=1}^m I_{T_j}$. In this case we can define a current $T'_j$ which can be viewed as the intersection of  $T_j$ and $T$ as follows.    Let $u_1, \ldots, u_m$ be local potentials of $T_1, \ldots, T_m$ respectively.  Let $T'_j:= \ddc (u_j|_V)$ for $1 \le j \le m$. One can see that $T'_j$ is independent of the choice of $u_j$, hence, $T'_j$ is a well-defined closed positive $(1,1)$-current on $V$. Let $\iota: V \hookrightarrow X$ be the natural inclusion. We can check that $\langle \wedge_{j=1}^m T_j \dot{\wedge} T \rangle$ is equal to the pushforward by $\iota$ of $\langle \wedge_{j=1}^m T'_j \rangle$. 
\end{remark}

\section{Monotonicity} \label{sec-mono}

In this section, we present a monotonicity property for relative non-pluripolar products. We begin with the following simple lemma.

\begin{lemma}\label{le-liminfnonpluripolar} Let $U$ be an open subset of $\C^n$. Let $(u^l_{j})_l$ be a sequence of psh functions on $U$ and $u_j$ a  psh function  on $U$ for $1 \le j \le m$. Let $(T_l)_l$  be a sequence of closed positive currents satisfying Condition $(*)$ and $T_l$ converges to a closed positive current $T$ on $U$ as $l \to \infty$.  Assume that one of the following conditions is satisfied:  

$(i)$  $u_{j}^l \ge u_j$ and $u_{j}^l$ converges to $u_j$ in $L^1_{loc}$ as $l \to\infty$. 

$(ii)$ $u_{j}^l$  increases to $u_j$ as $l \to \infty$ almost everywhere and $T$ has no mass on pluripolar sets.

Then, for every smooth weakly positive form $\Phi$ with compact support in $U$,  we have 
$$\liminf_{l \to\infty} \int_U \langle \wedge_{j=1}^m \ddc u_{j}^l \dot{\wedge} T_l  \rangle \wedge \Phi  \ge  \int_U \langle \wedge_{j=1}^m \ddc u_j \dot{\wedge} T \rangle \wedge \Phi.$$
\end{lemma}

When $T$ is the current of integration along $X$, a related statement in the compact setting was given in \cite[Theorem 2.3]{Lu-Darvas-DiNezza-mono}.

\proof   Let $\Phi$ be a smooth weakly positive form  with compact support in $U$. Assume $(i)$ holds. Let $u^l_{j  k}:= \max\{u_{j}^l, -k\}$ which converges to  $u_{j k}:= \max\{u_j, -k\}$ in $L^1_{loc}$ as $l \to \infty$ and $u^l_{jk} \ge u_{jk}$. Put 
$$R^l:= \langle \wedge_{j=1}^m \ddc u_{j}^l \dot{\wedge} T_l  \rangle, \quad R^l_k:=  \wedge_{j=1}^m \ddc u_{jk}^l \wedge T.$$
Similarly, we define $R, R_k$ by using the formula of $R^l,R^l_k$ respectively with $u_j^l, u_{jk}^l$ replaced by $u_j, u_{jk}$.   Since $u_j^l \ge u_j$, we get
\begin{align}\label{ine-sosanhujlvoiujliminf}
\int_U \bold{1}_{\cap_{j=1}^m  \{u_{j}^l > -k \}} R^l_k \wedge \Phi \ge \int_U \bold{1}_{\cap_{j=1}^m  \{u_{j} > -k \}} R_k^l \wedge \Phi.
\end{align}
By Theorem \ref{th-increasingsequenceMa},  we get $R^l_k \to R_k$ as $l \to \infty$. Using this together with the fact that  when $k$ is fixed, $u_{jk}^l$ is uniformly bounded in $l$, we see that  the strong uniform quasi-continuity for $u_{jk}$ with respect to $(T_l)_l$ (see Theorem \ref{th-capconvergedecreasing})  implies
\begin{align*} 
\liminf_{l \to \infty} \int_U \bold{1}_{\cap_{j=1}^m  \{u_{j} > -k \}} R^l_k \wedge \Phi &=\liminf_{l \to \infty} \int_U \bold{1}_{\cap_{j=1}^m  \{u_{jk} > -k \}} R^l_k \wedge \Phi \\
\nonumber
&  \ge  \int_U \bold{1}_{ \cap_{j=1}^m \{u_{jk} > -k \}} R_k \wedge \Phi=\int_U \bold{1}_{\cap_{j=1}^m  \{u_{j} > -k \}} R_k \wedge \Phi\\
\nonumber
&= \int_U \bold{1}_{\cap_{j=1}^m  \{u_{j} > -k \}} R \wedge \Phi.
\end{align*}
This combined with (\ref{ine-sosanhujlvoiujliminf}) yields
\begin{align*}
 \int_U \bold{1}_{\cap_{j=1}^m  \{u_{j} > -k \}} R \wedge \Phi &\le \liminf_{l \to \infty} \int_U \bold{1}_{\cap_{j=1}^m  \{u_{j}^l > -k \}} R^l_k \wedge \Phi \\
&\le \liminf_{l \to \infty} \int_U \bold{1}_{\cap_{j=1}^m  \{u^l_{j} > -k \}} R^l \wedge \Phi \le \liminf_{l \to \infty} \int_U R^l \wedge \Phi
\end{align*} 
for every $k$. Letting $k \to \infty$ in the last inequality and noticing that $R$ has no mass on $\cup_{j=1}^m \{u_j= -\infty\}$ give the desired assertion.

Now assume $(ii)$. Note that 
$$u^l_j \le u^{l+1}_j \le  u_j$$
 for every $l,j$. Thus, $\{u^l_j >-k\} \subset \{u^{l+1}_j >-k\} \subset \{u_j >-k\}$ and $\tilde{u}_j:= \lim_{l\to \infty} u^l_j \le u_j$. Using this and arguments similar to those in the previous paragraph,  we have
\begin{align*} 
\liminf_{l \to \infty} \int_U \bold{1}_{\cap_{j=1}^m  \{u^l_{j} > -k \}} R^l_k \wedge \Phi \ge  \int_U \bold{1}_{\cap_{j=1}^m  \{\tilde{u}_{j} > -k \}} R \wedge \Phi.
\end{align*}
Letting $k \to \infty$ gives
$$\liminf_{l \to \infty} \int_U R^l \wedge \Phi \ge  \int_U \bold{1}_{\cap_{j=1}^m  \{\tilde{u}_{j} > -\infty \}} R \wedge \Phi.$$
Recall that $\{u_j > \tilde{u}_j\}$ is pluripolar. This coupled with the fact that $T$ has no mass on pluripolar sets and Property $(iii)$ of  Proposition \ref{pro-sublinearnonpluripolar} yields
$$\liminf_{l \to \infty} \int_U R^l \wedge \Phi \ge  \int_U \bold{1}_{\cap_{j=1}^m  \{u_{j} > -\infty \}} R \wedge \Phi=\int_U R \wedge \Phi.$$
This finishes the proof. 
\endproof



Recall that for closed positive $(1,1)$-currents $R_1,R_2$ on $X$, we say that $R_1$ is \emph{less singular} than $R_2$ if for every local chart $U$ and psh function $w_j$ on $U$ such that   $R_j = \ddc w_j$ on $U$ for $j=1,2$, then $w_2 \le  w_1+ O(1)$ on compact subsets of $U$;  and $R_1, R_2$ are of \emph{the same singularity type} if $w_1 \le  w_2 + O(1)$ and $w_2 \le w_1+ O(1)$ on compact subsets of $U$.   The following generalizes  \cite[Proposition 2.1]{Lu-Darvas-DiNezza-mono}, see also \cite{WittNystrom-mono,Lu-cap-compare}. 

\begin{proposition} \label{pro-samesingularitytype}  Let $X$ be a compact complex manifold of dimension $n$. Let $m$ be an integer such that $1 \le m \le n$. Let $T$ be a closed positive current of bi-degree $(p,p)$ on $X$. Let $T_j, T'_j$ be closed positive $(1,1)$-currents on $X$ for $1 \le j \le m$ such that  $T_j, T'_j$ are of the same  singularity type and $T_j= \ddc u_j + \theta_j$, $T'_j: =\ddc u'_j+ \theta_j$, where $\theta_j$ is  a smooth form  and $u'_j, u_j$ are $\theta_j$-psh functions,  for every $1 \le j \le m$.   Assume that for every $J,J' \subset \{1, \ldots, m\}$ such that $J \cap J' = \varnothing$,  the product $\langle \wedge_{j \in J}T_{j} \wedge \wedge_{j'\in J'} T'_{j'} \dot{\wedge} T \rangle$ is well-defined. 
Then,  for every $\ddc$-closed smooth form $\Phi$, we have
\begin{align}\label{equa-samesingunonpluripolar}
\int_X \langle T'_{1} \wedge \cdots \wedge T'_{m} \dot{\wedge} T \rangle\wedge \Phi=  \int_X \langle T_{1} \wedge \cdots \wedge T_{m} \dot{\wedge} T \rangle\wedge \Phi.
\end{align}
\end{proposition} 

\proof  By compactness of $X$, we can assume $u_j,u'_j \le 0$. By the hypothesis, $\{u_j= -\infty\}= \{u'_j= -\infty\}$ and  $w_j:= u_j - u'_j$ is bounded outside $\{u_j= -\infty\}$.  Without loss of generality, we can assume that 
\begin{align}\label{ine-hieuujujphay}
|w_j| \le 1
\end{align}
outside $\{u_j= -\infty\}$.  Let $A:= \cup_{j=1}^m \{u_j= -\infty\}$ which is a complete pluripolar set. Put $u_{j k}:= \max\{u_j, -k\}$, $u'_{j k}:= \max\{u'_j, -k\}$ and
\begin{align}\label{equa-defvarphik}
\psi_k:= k^{-1}\max\{\sum_{j=1}^n (u_j+u'_j), -k\}+ 1
\end{align}
which is quasi-psh and $0\le \psi_k \le 1$, $\psi_k(x)$ increases to $1$ for $x\not \in A$. We have $\psi_k(x) =0$ if $u_j(x)\le -k$ or $u'_j(x) \le -k$ for some $j$. 
Put $w_{j k}:= u_{j k} - u'_{j k}$. By (\ref{ine-hieuujujphay}), we have
\begin{align}\label{ine-hieuujujphay2}
|w_{j k}| \le 1
\end{align}
on $X$.   Let $J, J' \subset \{1, \ldots, n\}$ with $J \cap J' = \varnothing$ and 
$$R_{JJ'k}:= \wedge_{j \in J}  (\ddc u_{j k} + \theta_j) \wedge \wedge_{j' \in J'} (\ddc u'_{j' k} + \theta_{j'})\wedge T.$$
The last current is the difference of two closed positive $(|J|+|J'|+p, |J|+|J'|+p)$-currents. Hence, $R_{JJ'k}$ might not be positive in general and it is not clear how to control its mass  as $k \to \infty$.  This is a subtle point which we need to pay attention to. The relative non-pluripolar product
$$R_{JJ'}:= \big \langle \wedge_{j \in J}  (\ddc u_{j} + \theta_j) \wedge \wedge_{j' \in J'} (\ddc u'_{j'} + \theta_{j'}) \dot{\wedge}T\big\rangle$$
 exists by our assumption.  Let 
$$B_k:= \cap_{j \in J} \{u_j > -k\} \cap \cap_{j' \in J'}\{ u'_{j'} > -k\}.$$
 By Lemma \ref{le-globalpotentialnonpluri}, we get
$$0\le \bold{1}_{B_{k}} R_{JJ'}= \bold{1}_{B_{k}} R_{J J' k}$$
for every $J,J',k$.  Put $\tilde{R}_{JJ'}:= \bold{1}_{X \backslash A} R_{JJ'}$. The last current is closed (for example see \cite[Remark 1.10]{BEGZ}) and positive because $R_{JJ'} \ge 0$.  Using  the fact that $\{\psi_{k} \not =0\} \subset B_{k} \backslash A$, we get 
\begin{align}\label{eq-bieudienRjkquaRjnew}
\psi_{k} \tilde{R}_{JJ'}= \psi_{k} R_{JJ'}= \psi_{k} R_{JJ'k}.
\end{align}
 
\noindent
\textbf{Claim. } Let $j'' \in \{1, \ldots, m\} \backslash (J \cup J')$. Let  $\Phi$ be  $\ddc$-closed smooth form of  bi-degree $(p',p')$ on $X$, where $p':=n- |J|-|J'|-p-1$.  Then 
\begin{align}\label{limit-wjkbichanbansua}
\lim_{k \to \infty}\int_X \psi_{k} \ddc w_{j'' k} \wedge R_{JJ'k} \wedge \Phi =0.
\end{align}

\noindent
We prove Claim. Let $\omega$ be a Hermitian metric on $X$. 
Let $\eta: =\sum_{j=1}^m 2\, \theta_j$. We have  $\ddc\psi_k+ k^{-1}\eta \ge 0$ for every $k$. By integration by parts and (\ref{eq-bieudienRjkquaRjnew}),
\begin{align}\label{eq-tachtichphanRjjkkweeeeee}
\int_X\psi_{k}\ddc w_{j'' k} \wedge R_{JJ'k}\wedge \Phi=\int_X w_{j'' k} \ddc (\psi_{k}R_{JJ'k}\wedge \Phi)=\int_X w_{j'' k} \ddc (\psi_{k} \tilde{R}_{JJ'}\wedge \Phi).
\end{align}
Observe that 
$$\ddc (\psi_{k}  \tilde{R}_{JJ'}\wedge \Phi)= \ddc\psi_k \wedge \Phi \wedge  \tilde{R}_{JJ'}+  2d\psi_k \wedge  \dc \Phi \wedge  \tilde{R}_{JJ'}$$
because $\ddc \Phi = 0$ and $\Phi$ is of bi-degree $(p',p')$. Write $\dc \Phi$ locally as a complex linear combination of forms like $\tau_j \wedge \Phi_j$, where $\tau_j$ is a $(0,1)$-form or a $(1,0)$-form  and $\Phi_j$ is  a positive form. Hence, we can use the  Cauchy-Schwarz inequality to obtain
\begin{align*}
\bigg| \int_X w_{j'' k}\, d\psi_k \wedge  \dc \Phi \wedge  \tilde{R}_{JJ'} \bigg|&\le  \bigg(\int_X |w_{j'' k}^2| d\psi_k \wedge  \dc\psi_k \wedge  \tilde{R}_{JJ'}\wedge \Phi_0\bigg)^{\frac{1}{2}}  \bigg(\int_X   \tilde{R}_{JJ'}\wedge \Phi_0 \wedge \omega\bigg)^{\frac{1}{2}},
\end{align*}
where $\Phi_0:= \omega^{n- |J|- |J'|-1}$. We deduce that 
$$\bigg| \int_X w_{j'' k}\, d\psi_k \wedge  \dc \Phi \wedge  \tilde{R}_{JJ'} \bigg|   \lesssim  \| \tilde{R}_{JJ'}\|^{\frac{1}{2}}\bigg(\int_X d\psi_k \wedge  \dc\psi_k \wedge  \tilde{R}_{JJ'}\wedge \Phi_0\bigg)^{\frac{1}{2}}$$
by (\ref{ine-hieuujujphay2}).  Recall that  $\{\lim_{k \to\infty} \psi_k <1\}$ is equal to the complete pluripolar set $A$. Using this, Remark \ref{re-th-increasingsequenceMa} and the fact that $\tilde{R}_{JJ'}$ has no mass on $A$,  we get
$$\lim_{k \to \infty} d\psi_k \wedge  \dc\psi_k \wedge  \tilde{R}_{JJ'}=\lim_{k \to \infty} (\ddc \psi^2_k - \psi_k \ddc \psi_k) \wedge  \tilde{R}_{JJ'}=0$$
(we recall  that to get  Remark \ref{re-th-increasingsequenceMa} for $\tilde{R}_{JJ'}$ and $\psi_k$,  we need to use the strong quasi-continuity of $\psi_k$ with respect to the capacity associated to $\tilde{R}_{JJ'}$).  Thus we obtain
\begin{align}\label{limit-wjkbichan00}
\lim_{k \to \infty }\int_X w_{j'' k}\,  d\psi_k \wedge  \dc \Phi \wedge  \tilde{R}_{JJ'} =0.
\end{align}
On the other hand,  since 
$$ \int_X w_{j'' k} \ddc\psi_k \wedge \Phi \wedge  \tilde{R}_{JJ'}=  - \int_X d w_{j'' k} \wedge \dc\psi_k \wedge \Phi \wedge  \tilde{R}_{JJ'}+  \int_X w_{j'' k} \dc\psi_k \wedge d \Phi \wedge  \tilde{R}_{JJ'},$$
 using similar arguments, we get 
$$\lim_{k' \to \infty}\bigg| \int_X w_{j'' k} \ddc\psi_k \wedge \Phi \wedge  \tilde{R}_{JJ'}\bigg|=0 $$
Combining this with (\ref{limit-wjkbichan00}) and (\ref{eq-tachtichphanRjjkkweeeeee}) yields  (\ref{limit-wjkbichanbansua}). Claim follows.

Now let $S:= \langle T_{1} \wedge \cdots T_{n}\dot{\wedge} T \rangle - \langle T'_{1} \wedge \cdots T'_{n}\dot{\wedge} T \rangle$. Using $T_{jk}= T'_{jk}+ \ddc w_{jk}$, one can check that
\begin{align*}
\int_X \psi_{k} S\wedge \Phi &= \int_X \psi_{k} \wedge_{j=1}^m T_{jk}\wedge T\wedge \Phi-\int_X \psi_{k} \wedge_{j=1}^m T'_{jk}\wedge T\wedge \Phi\\
&=\sum_{s=1}^m \int_X \psi_{k} \wedge_{j=1}^{s-1} T'_{jk}\wedge \ddc w_{s k} \wedge \wedge_{j=s+1}^m T_{jk} \wedge T \wedge \Phi.
\end{align*}
This together with Claim yields  $\langle S, \Phi \rangle = \lim_{k\to \infty} \langle \psi_{k}S, \Phi\rangle =0$
 for every $\ddc$-closed smooth  $\Phi$.  This finishes the proof.
\endproof

\begin{remark}  Our proof of Proposition \ref{pro-samesingularitytype} still works if  $T'_j, T_j$ are not in the same cohomology class. In this case, one just needs to modify  (\ref{equa-samesingunonpluripolar})  accordingly.
\end{remark}


The following result is the monotonicity property of relative non-pluripolar products mentioned in Introduction.

\begin{theorem} \label{th-monoticity} Let $X$ be a compact K\"ahler manifold of dimension $n$.  Let $T_j, T'_j$ be closed positive $(1,1)$-currents on  $X$ for $1 \le j \le m $ such that $T_j, T'_j$ are in the same cohomology class for every $j$ and $T'_j$ is less singular than $T_j$ for $1 \le j \le m$.  Let $T$ be a closed positive current on $X$.   Then, we have 
\begin{align}\label{ine-lesssingu-ddcclaosenonpluri}
\big\{  \langle T_1 \wedge \cdots \wedge T_m  \dot{\wedge} T\rangle \big\} \le \big\{ \langle T'_1 \wedge \cdots \wedge T'_m \dot{\wedge} T \rangle  \big\}.
\end{align} 
\end{theorem}


\proof  Write $T_j= \ddc u_j + \theta_j$, $T'_j= \ddc u'_j + \theta_j$. Without loss of generality, we can assume that $u'_j \ge u_j$.  For $l \in \N$, put $u_{j}^l:= \max\{ u_j, u'_j - l \}$ which is of the same singularity type as $u'_j$. Notice that $\ddc u_{j}^l+ \theta_j \ge 0$. Since $X$ is K\"ahler, the current  $\langle \wedge_{j=1}^m (\ddc u_{j}^l+ \theta_j)  \dot{\wedge} T\rangle$  is of mass uniformly bounded in $l$.  Let $S$ be a limit current of $\langle \wedge_{j=1}^m (\ddc u_{j}^l+ \theta_j)  \dot{\wedge} T\rangle$ as $l \to \infty$.  Since $u_{j}^l$ decreases to $u_j$ as $l \to \infty$, we can apply Lemma \ref{le-liminfnonpluripolar} to get
$$S \ge \langle \wedge_{j=1}^m T_j \dot{\wedge} T \rangle.$$
Consequently, $\{S\} \ge  \{\langle \wedge_{j=1}^m T_j \dot{\wedge} T \rangle\}$.  Using Proposition \ref{pro-samesingularitytype}, we see that $\{S\}$ is equal to $\{\langle \wedge_{j=1}^m T'_j  \dot{\wedge} T \rangle\}$. The desired assertion, hence,  follows. This finishes the proof. 
\endproof

\begin{remark} \label{re-Tofbidegre11mono} Let the notation be as in Theorem \ref{th-monoticity}.  Let $T$ be of bi-degree $(1,1)$ and $T'$ a closed positive $(1,1)$-current which is less singular than $T$. Then, by using arguments similar to those in the proof of Theorem \ref{th-monoticity}, we  can prove that
\begin{align}\label{ine-lesssingu-ddcclaosenonpluribac11}
\big\{  \langle T_1 \wedge \cdots \wedge T_m  \dot{\wedge} T\rangle \big\} \le \big\{ \langle T'_1 \wedge \cdots \wedge T'_m \dot{\wedge} T' \rangle  \big\}
\end{align} 
(here we need to use Lemma \ref{le-liminfnonpluripolar} for a suitable sequence $(T_l)_l$ provided by Theorem \ref{th-decreasing-unbounded}). The inequality (\ref{ine-lesssingu-ddcclaosenonpluribac11}) offers us a way to define a notion of full mass intersection when $T$ is of bi-degree $(1,1)$. This notion differs from those used below and in \cite{BEGZ}, albeit all of them are closely related. We will not go into details in this paper.
\end{remark}

Consider, from now on,  a compact K\"ahler manifold $X$ with a K\"ahler form $\omega$. Let $T$ be a closed positive $(p,p)$-current on $X$.  For every pseudoeffective $(1,1)$-class $\beta$ in $X$, we define its \emph{polar locus $I_\beta$} to be that of a current with minimal singularities in $\beta$. This is independent of the choice of a current with minimal singularities. 

  Let $\alpha_1, \ldots,\alpha_m$ be pseudoeffective $(1,1)$-classes of $X$. Let $T_{1, \min}, \ldots, T_{m,\min}$ be currents with minimal singularities in  the classes $\alpha_1, \ldots, \alpha_m$ respectively.   By Theorem \ref{th-monoticity} and Lemma \ref{le-compactKahlernonpluripolar}, the class $\{\langle T_{1,\min} \wedge \cdots \wedge T_{m, \min} \dot{\wedge} T \rangle \}$ is  a well-defined pseudoeffective class which is independent of the choice of $T_{j, \min}$. We denote the last class by $\{\langle \alpha_1\wedge  \cdots \wedge \alpha_m \dot{\wedge} T \rangle\}$.  

\begin{proposition} \label{pro-productclasss}
$(i)$ The product  $\{\langle  \wedge_{j=1}^m \alpha_j \dot{\wedge} T \rangle\}$ is symmetric and   homogeneous  in $\alpha_1,$ $\ldots, \alpha_m$.  

$(ii)$ Let  $\alpha'_1$ are a pseudoeffective $(1,1)$-class.  Assume that $T$ has no mass on $I_{\alpha_1}\cup I_{\alpha'_1}$.  Then, we have 
$$\{\langle (\alpha_1+ \alpha'_1) \wedge \wedge_{j=2}^m \alpha_j \dot{\wedge} T \rangle\} \ge \{\langle  \wedge_{j=1}^m \alpha_j\wedge T \rangle\}+ \{\langle  \alpha'_1 \wedge \wedge_{j=2}^m \alpha_j \dot{\wedge} T \rangle\}.$$

$(iii)$ Let $1 \le l \le m$ be an integer. Let  $\alpha''_1, \ldots, \alpha''_l$ be a pseudoeffective $(1,1)$-class such that $\alpha''_j \ge \alpha_j$ for $1 \le j \le l$.  Assume that $T$ has no mass on $I_{\alpha''_j- \alpha_j}$ for every $1 \le j \le l$. Then, we have 
$$\{\langle \wedge_{j=1}^l \alpha''_j \wedge  \wedge_{j=l+1}^m \alpha_j \dot{\wedge} T \rangle\} \ge  \{ \langle \wedge_{j=1}^m \alpha_j \dot{\wedge} T \rangle\}.$$

$(iv)$ If $T$ has no mass on proper analytic subsets on $X$, then the product   $\{\langle  \wedge_{j=1}^m \alpha_j \dot{\wedge} T \rangle\}$ is continuous on the set of $(\alpha_1,\ldots, \alpha_m)$ such that $\alpha_1, \ldots, \alpha_m$ are big.

$(v)$  If $T$ has no mass on proper analytic subsets on $X$ and  $\alpha_1, \ldots, \alpha_m$ are big nef, then we have    
$$\{\langle \wedge_{j=1}^m \alpha_j \dot{\wedge} T \rangle\}= \wedge_{j=1}^m \alpha_j \wedge \{T\}.$$

\end{proposition}

We refer to \cite{Boucksom-derivative-volume,BEGZ} for related statements in the case where $T$ is the current of integration along $X$. 

\proof  The  desired assertion $(i)$ follows from Proposition \ref{pro-sublinearnonpluripolar}. 
We now prove $(ii)$.  Let $T_{\min, \alpha_j}, T_{\min, \alpha'_j}$ be currents with minimal singularities in $\alpha_j, \alpha'_j$ respectively. Observe that $T_{\min,\alpha_j}+ T_{\min, \alpha'_j}$ are in $(\alpha_j+ \alpha'_j)$. Thus, by Theorem \ref{th-monoticity}, we get 
$$\{\langle  (\alpha_1+ \alpha'_1) \wedge  \wedge_{j=2}^m \alpha_j  \dot{\wedge} T \rangle\}\ge \{\langle (T_{\min,\alpha_1}+ T_{\min, \alpha'_1}) \wedge \wedge_{j=2}^m T_{\min,\alpha_j}\dot{\wedge}T \rangle \}.$$
 The last class is equal to $\{\langle  \wedge_{j=1}^m T_{\min,\alpha_j}\dot{\wedge} T \rangle\}+ \{\langle  T_{\min, \alpha'_1} \wedge \wedge_{j=2}^m T_{\min,\alpha_j} \dot{\wedge} T \rangle\}$ because of the hypothesis and Property $(iv)$ of Proposition \ref{pro-sublinearnonpluripolar}. Hence, $(ii)$ follows. Similarly, we get $(iii)$ by using Property $(v)$ of Proposition \ref{pro-sublinearnonpluripolar}.

We prove $(iv)$.   Observe that by a result of Demailly on analytic approximation of currents (\cite{Demailly_analyticmethod}), the polar locus $I_\beta$ of a big class $\beta$ is contained in a proper analytic subset of $X$ if $\alpha$ is big. Using this, we see that $(iv)$ is a direct consequence of $(iii)$ and the observation that given a constant $\epsilon>0$,  for every $\alpha'_j$ closed enough to $\alpha_j$, we have  that the classes $\alpha_j' -(1-\epsilon) \alpha_j$ and $ (1+ \epsilon) \alpha_j- \alpha'_j$ are big (we use here the bigness of $\alpha_j$).  

It remains to check $(v)$. By $(iv)$ and the bigness of $\alpha_j$, we get 
$$\lim_{\epsilon \to 0} \{\langle \wedge_{j=1}^m (\alpha_j+ \epsilon \{\omega\}) \dot{\wedge} T\rangle \} = \{\langle \wedge_{j=1}^m \alpha_j \dot{\wedge} T \rangle\}.$$
Since $\alpha_j$ is nef, the limit in the left-hand side of the last equality is equal to $\wedge_{j=1}^m \alpha_j \wedge \{T\}$.  The proof is finished.
\endproof


When $T$ is the current of integration along $X$, we write $\langle \alpha_1 \wedge \cdots \wedge \alpha_m \rangle$ for $\{\langle \alpha_1\wedge  \cdots \wedge \alpha_m \dot{\wedge} T \rangle\}$.  
We would like to comment that  the class $\langle \alpha_1\wedge  \cdots \wedge \alpha_m \rangle$ is always bounded from above by the \emph{positive product } of $\alpha_1, \ldots, \alpha_m$ defined in  \cite[Definition 1.17]{BEGZ}. They are equal  if $\alpha_1, \ldots, \alpha_m$ are big by Property $(iv)$ of Proposition \ref{pro-productclasss}. However, we don't know if they are equal in general, even if $\alpha_1, \ldots, \alpha_m$ are nef.  

\begin{question} Given nef classes $\alpha_1, \ldots, \alpha_m$, is $\langle \alpha_1 \wedge \cdots \wedge \alpha_m \rangle$ defined above equal to the positive product of $\alpha_1, \ldots, \alpha_m$ introduced in \cite{BEGZ}? 
\end{question}

 Let $T_1, \ldots, T_m$ be closed positive $(1,1)$-currents on $X$.  By Theorem \ref{th-monoticity},  we  have $\{\langle \wedge_{j=1}^m T_j \dot{\wedge} T \rangle  \} \le \{\langle \wedge_{j=1}^m \{T_j \} \dot{\wedge} T\rangle \}$. The equality occurs if the masses of these two classes are equal. This is the reason for the following definition. 

\begin{definition} \label{def-fullmassinter}  We say that $T_1, \ldots, T_m$ are of \emph{full mass (non-pluripolar)  intersection relative to $T$} if we have $\{\langle \wedge_{j=1}^m T_j \dot{\wedge} T \rangle  \} = \{\langle \wedge_{j=1}^m \{T_j \} \dot{\wedge} T\rangle \}$. 
\end{definition}

When  $T$ is the current of integration along $X$, we simply say  ``full mass intersection" instead of ``full mass intersection relative to $T$". In the last case, we underline that this notion is  the one given in \cite{BEGZ} if $\{T_1\}, \ldots, \{T_m\}$ are big. In general, if $T_1, \ldots, T_m$ are of full mass intersection in the sense of \cite{BEGZ}, then they are so in the sense of Definition \ref{def-fullmassinter}. However, we don't know whether the reversed statement holds. 

  Let $\mathcal{E}(\alpha_1,\ldots, \alpha_m,T)$ be the set of $(T_1, \ldots, T_m)$ such that $T_j \in \alpha_j$ for $1 \le j \le m$ and  
$$\{\langle \wedge_{j=1}^m T_j \dot{\wedge} T \rangle  \} = \{\langle  \wedge_{j=1}^m \alpha_j \dot{\wedge} T\rangle \},$$
 or equivalently  
$$\int_X \langle \wedge_{j=1}^m T_j \dot{\wedge} T \rangle \wedge \omega^{n-m-p} = \int_X \langle \wedge_{j=1}^m \alpha_j \dot{\wedge} T \rangle \wedge \omega^{n-m-p}.$$ 
Several less general versions of full mass intersections were introduced in \cite{GZ-weighted,BEGZ,Lu-Darvas-DiNezza-mono}. Note that if $T_j= \ddc u_j + \theta_j$ for some smooth K\"ahler form $\theta_j$ and $u_j$ a $\theta_j$-psh function for every $j$, then $T_1,\ldots, T_m$ are of full mass intersection relative to $T$ if and only if 
\begin{align} \label{ine-dkfullmasscutoffk}
\int_{\cup_{j=1}^m\{ u_j \le -k\}} \wedge_{j=1}^m T_{jk} \wedge T\wedge \omega^{n-m-p} \to 0 
\end{align}
as $k \to \infty$, where $T_{jk}:= \ddc \max\{u_j, -k\}+ \theta_j$ for $1 \le j \le m$.

 The following result tells us that currents with full mass intersection satisfy the convergence along decreasing or increasing sequences of potentials.  

\begin{theorem} \label{th-convergencMAgenernonpluri}  Let $X$ be a compact K\"ahler  manifold of dimension $n$.
Let $T_{jl}$ be  a closed positive $(1,1)$-current for $1 \le j \le m,$ $l \in \N$ such that $T_{jl}$ converges to $T_j$ as $l \to \infty$. Let $(U_s)_s$ be a finite covering of $X$ by open subsets such that $T_{jl} = \ddc u_{jl,s}$, $T_j= \ddc u_{j,s}$ on $U_s$ for every $s$. Assume that 

$(i)$ $T_1, \ldots, T_m$ are of full mass intersection relative to $T$,

$(ii)$  $\{\langle  \wedge_{j=1}^m\alpha_{jl}\dot{\wedge} T \rangle\}$  converges to $\{\langle \wedge_{j=1}^m \alpha_j  \dot{\wedge} T \rangle\}$ as $l \to \infty$, where $\alpha_{jl}:= \{T_{jl}\}$, $\alpha:= \{T_j\}$

and one of the following two conditions hold:

$(iii)$ $u_{jl,s} \ge u_{j,s}$ and $u_{jl,s}$ converges to $u_{j,s}$ in $L^1_{loc}$ as $l \to\infty$, 

$(iii')$ $u_{jl,s}$  increases to $u_{j,s}$ as $l \to \infty$ almost everywhere and $T$ has no mass on pluripolar sets. 

Then we have 
$$\langle T_{1l} \wedge \cdots \wedge T_{ml}  \dot{\wedge}  T \rangle \to \langle T_1 \wedge \cdots \wedge T_m  \dot{\wedge}  T \rangle$$
as $l \to \infty$.
\end{theorem}

\proof Observe that $\langle \wedge_{j=1}^m T_{jl} \dot{\wedge} T \rangle $ is of uniformly bounded mass in $l$ by $(ii)$.   Let $S$ be a limit current of $\langle \wedge_{j=1}^m T_{jl} \dot{\wedge} T \rangle $ as $l \to \infty$. By $(ii)$ again, we get $\{S\} \le \{\langle \wedge_{j=1}^m \alpha_j  \dot{\wedge} T \rangle\}$.  By Lemma \ref{le-liminfnonpluripolar} and Condition $(iii)$ or $(iii')$, we have 
$$ S \ge  \langle \wedge_{j=1}^m T_j \dot{\wedge} T \rangle.$$
By $(i)$, the current in the right-hand side is of mass equal to $\langle \wedge_{j=1}^m \alpha_j \dot{\wedge} T \rangle$ which is greater than or equal to the mass of  $S$. So we get the equality $S=  \langle \wedge_{j=1}^m T_{j} \dot{\wedge} T \rangle$. The proof is finished.
\endproof

Note that we don't require that  $T_{1l}, \ldots, T_{ml}$ are of full mass intersection relative to $T$. However, if $T_{jl}, T_j$ are in the same cohomology class for every $j,l$, then by Theorem \ref{th-monoticity}, Conditions $(iii)$ and $(i)$ of Theorem \ref{th-convergencMAgenernonpluri}  imply  that $T_{1l}, \ldots, T_{ml}$ are of full mass intersection relative to $T$  for every $l$. We notice here certain similarity of this result with \cite[Propsition 5.4]{DinhSibony_pullback}, where a notion of full mass intersection was considered for intersections of currents with analytic sets.   

Here are some basic properties of currents with relative full mass intersection. 

\begin{lemma} \label{le-hasnomassonITjkahler} Let $X$ be a compact K\"ahler manifold. Let $T_1, \ldots, T_m,T$ be closed positive currents on $X$ such that $T_j$ is of bi-degree $(1,1)$ and the cohomology class of $T_j$ is K\"ahler for every $j$. Then,  if $T_1, \ldots, T_m$ are of full mass intersection relative to $T$, then the following three properties hold.

$(i)$  $T$ has no mass on $\cup_{j=1}^m I_{T_j}$.

$(ii)$ Let $T'_j$ be closed positive $(1,1)$-currents whose cohomology class are K\"ahler for $1 \le j \le m$ such that $T'_j$ is less singular than $T_j$ for $1 \le j \le m$. Then  $T'_1, \ldots, T'_m$ are of full mass intersection relative to $T$.

$(iii)$ For every subset $J=\{j_1, \ldots, j_{m'}\} \subset \{1, \ldots, m\}$,  the currents $T_{j_1}, \ldots, T_{j_{m'}}$ are of full mass intersection relative to $T$. Moreover, if $T$ is the current of integration along $X$, then  $T_j$ has no mass on $I_{T_j}$ for $1 \le j \le m$. 

Moreover, we have 

$(iv)$  $T_1, \ldots, T_m$ are of full mass intersection relative to $T$ if and only if $T_1+C \omega, \ldots, T_m+C \omega$ are of full mass intersection relative to $T$  for every constant $C>0$.
\end{lemma}

\proof The desired property $(i)$ is a direct consequence of $(vii)$ of Proposition \ref{pro-sublinearnonpluripolar} and the fact that $\{T_j\}$ is K\"ahler for every $j$. 
Note that if $\{T'_j\}=\{T_j\}$ for every $j$, then $(ii)$ is a direct consequence of  the monotonicity of relative non-pluripolar products (Theorem \ref{th-monoticity}). The first claim of the desired property $(iii)$ is deduced from the last assertion applies to $T'_j:= T_j$ for $j \in J$ and $T'_j:= \theta_j$ for $j \not \in J$, where  $\theta_j$ is a smooth K\"ahler form in the class $\{T_j\}$. In particular,  we have $\{\langle T_j \rangle\} = \{T_j\}$. Recall that $\langle T_j \rangle =\bold{1}_{X \backslash I_{T_j}} T_j$ (Proposition \ref{pro-tinhnonpluripolarclassically}). Hence, 
$$\| \langle T_j \rangle\|+ \| \bold{1}_{I_{T_j}} T_j\|  = \|T_j \|.$$
It follows that $\bold{1}_{I_{T_j}} T_j=0$ or equivalently, $T_j$ has no mass on $I_{T_j}$. Hence $(iii)$ follows.  

Now observe that $(iv)$ is a direct consequence of $(iii)$, Property $(iv)$ of Proposition \ref{pro-sublinearnonpluripolar} and $(i)$. Finally, when $\{T'_j\} \not = \{T_j\}$, Property $(iv)$ allows us to  add to $T_j,T'_j$ suitable K\"ahler forms such that they are in the same cohomology class. So $(ii)$ follows. The proof is finished.
\endproof

In the last part of this section, we study Lelong numbers of currents of relative full mass intersection.

\begin{lemma} \label{le-reduction-analyticsingularity} Let $u$ be a quasi-psh function on $X$.  Let $V$ be an analytic subset of $X$.  Let $v$ be a quasi-psh function on $X$ with analytic singularities along  $V$. Assume that for every $x \in V$, we have $\nu(u,x)> 0$.  Then there exists a constant $c >0$ such that $u \le  c v+ O(1)$ on $X$.
\end{lemma}

Recall that $v$ is said to \emph{have analytic singularities along $V$} if locally $v= c \log \sum_{j=1}^l |f_j|+ w$, where $c>0$ is a constant, $f_j$ is holomorphic for every $j$ such that $V$ is locally equal to  $\{f_j=0: 1 \le j \le l\}$, and  $w$ is a bounded Borel function, see \cite[Definition 1.10]{Demailly_analyticmethod}. Moreover, given an analytic set $V$ in $X$, there always exists a quasi-psh  function $v$ having analytic singularities along $V$.  We can construct such a function by using a partition of unity and local generators defining $V$, see \cite[Lemma 2.1]{Demailly-Paun}.

\proof We can assume that $u,v$ are $\omega$-psh functions.  Put $T:= \ddc u+\omega$. Note that by hypothesis and Siu's semi-continuity theorem,  there is a constant $c_0>0$ such that  $\nu(u,x)\ge c_0>0$ for every $x \in V$ and $\nu(u,x)= c_0$ for generic $x$ in $V$. Consider first the case where $V$ is of codimension 1. Thus $\bold{1}_{V}T= c_0 [V]$ is a non-zero positive current (\cite[Lemma 2.17]{Demailly_analyticmethod}) and $c_1 v_1 + O(1) \le v \le c_1 v_1+ O(1)$  for some constant $c_1>0$, where $v_1$ is a potential of $[V]$. Since $T \ge \bold{1}_{V}T= c_0[V]$, we get $u \le  c_0 v_1+O(1)$ on $X$. Thus $u \le  (c_0/c_1) \, v + O(1)$ on $X$.  

Consider now $\codim V  \ge 2$.   By desingularizing $V$, we obtain a compact K\"ahler manifold $X'$ and a surjective map $\rho: X' \to X$  such that $V':=\rho^{-1}(V)$ is a  hypersurface.  Note that $\nu(u \circ \rho, x)>0$ for every $x \in V'$ and $v\circ \rho$ is a current with analytic singularities along $V'$. Applying the first part of the proof gives the desired assertion.  This finishes the proof.
\endproof

Recall for every psh function $v$,  the unbounded locus $L(v)$ of $v$ defined to be the set of $x$ such that $v$ is unbounded in every open neighborhood of $x$. Observe that if $v$ has analytic singularities along $V$, then $L(v)=V$.  For every closed positive $(1,1)$-current $T$, we define $L(T)$ to be the unbounded locus of a potential of $T$. Here is a necessary condition for  currents to be of relative full mass intersection in terms of Lelong numbers. 

\begin{theorem} \label{the-lelongobstructiontongmpnhonhonn}
Let $X$ be a compact K\"ahler  manifold. Let $T_1, \ldots, T_m$ be closed positive $(1,1)$-currents on $X$ such that the cohomology class of $T_j$ is K\"ahler for $1 \le j \le m$  and $T$ a closed positive $(p,p)$-current on $X$ with $p+m \le n$.   Let $V$ be an irreducible analytic subset in $X$ such that the generic Lelong numbers of $T_1, \ldots, T_m,T$ are strictly positive. Assume $T_1, \ldots, T_m$ are of full mass intersection relative to $T$. Then, we have  $\dim V <n-p-m$. 
\end{theorem}

\proof 
 Since $\theta_j$ is K\"ahler for every $j$, we can use Theorem \ref{th-monoticity}, Lemma \ref{le-reduction-analyticsingularity} and the comment following it to reduce the setting to the case where $T_j$ are currents with analytic singularities along $V$, for $1 \le j \le m$ (hence $L(T_j)= V$).  Since $T$ is of bi-degree $(p,p)$ and $T$ has positive Lelong number everywhere on $V$, we deduce that the dimension of $V$ is at most $n-p$. Let $s$ be a nonnegative integer such that $\dim V =  n-p -s$. We need to prove that $s >m$.  Suppose on the contrary that $s \le m$. By 
Lemma \ref{le-hasnomassonITjkahler}, the currents $T_1,\ldots, T_s$ are of full mass intersection relative to $T$. 

Since $L(T_j)= V$ for every $j$, we see that for every $J \subset \{1, \ldots, s\}$, the Hausdorff dimension of  $\cap_{j \in J}L(T_j) \cap \supp T$ is less than or equal to that of $\cap_{j \in J}L(T_j)$ which is equal to $\dim V= n-p- s \le n-p - |J|$. Thus, by \cite{Demailly_ag,Fornaess_Sibony}, the intersection of $T_1, \ldots, T_s, T$ is classically well-defined.  By a comparison result on Lelong numbers (\cite[Page  169]{Demailly_ag}) and the fact that $T_1, \ldots, T_m,T$ have strictly positive Lelong number at every point in $V$, we see that  the Lelong number of $\wedge_{j=1}^sT_j \wedge T$ at every point of $V$ is strictly positive. Thus, the current $\wedge_{j=1}^sT_j \wedge T$ has strictly positive mass on $V$ (see  \cite[Lemma 2.17]{Demailly_analyticmethod}). So by Proposition \ref{pro-tinhnonpluripolarclassically},  $\langle \wedge_{j=1}^sT_j \wedge T \rangle$ is not  of maximal mass. This is a contradiction. 
 This finishes the proof.
\endproof

\section{Weighted class of currents of  relative full mass intersection} \label{sec-weightedclass}

In this subsection, we introduce the notion of weighted classes of currents with  full mass intersection relative to a closed positive current $T$. We only consider convex weights in the sequel.  We refer to  \cite{GZ-weighted,BEGZ,Darvas_book} for the case where $T \equiv 1$ and informations on  non-convex weights. We also note that \cite{BEGZ} considers currents in big classes whereas our setting here is restricted to currents in K\"ahler classes. 

We fix our setting. Let $\omega$ be a K\"ahler form on $X$. Let $m \in \N^*$.   Let $T_j= \ddc u_j + \theta_j$ be a closed positive $(1,1)$-current where $\theta_j$ is a \emph{smooth K\"ahler} form and $u_j$ is a negative $\theta_j$-psh function for $1 \le j \le m$.   Let $0 \le p \le n$ be an integer and $T$ a closed positive current of bi-degree $(p,p)$ on $X$.  We assume $p+m \le n$. This is a minimal assumption because otherwise the relative non-pluripolar product is automatically zero by a bi-degree reason. 

For $k \in \N$, put $u_{jk}:= \max\{u_j, -k\}$ and $T_{jk}:= \ddc u_{jk}+ \theta_j$. Note that $T_{jk} \ge 0$ because $\theta_j \ge 0$.   Let $\beta$ be a K\"ahler $(1,1)$-class and $\mathcal{E}_m(\beta,T)$  the set of closed positive $(1,1)$-currents $P$ in the class $\beta$ such that $(P, \ldots, P)$ ($m$ times $P$) is in $ \mathcal{E}(\beta, \ldots, \beta, T)$ ($m$ times $\beta$).  The class $\mathcal{E}_n(\beta)$ was introduced in \cite{GZ-weighted}.   We first prove the following result giving the convexity of the class of currents of relative full mass intersection.

\begin{theorem} \label{the-mathcalEnalpha}  Assume that  $T_j, \ldots, T_j$  ($m$ times $T_j$) are of full mass intersection relative to $T$  for $1 \le j \le m$. Then,  $T_1, \ldots, T_m$ are also of full mass intersection relative to $T$. In particular,  for  K\"ahler  $(1,1)$-classes $\beta, \beta'$, we have
 $$\mathcal{E}_{m}(\beta, T)+ \mathcal{E}_{m}(\beta',T) \subset \mathcal{E}_{m}(\beta+ \beta', T)$$
and $\mathcal{E}_{m}(\beta,T)$ is convex . 
\end{theorem}

\proof Note that the second desired assertion is a direct consequence of the first one (recall that $\{T_j\}$'s are K\"ahler). By Lemma \ref{le-hasnomassonITjkahler}, $T$ has no mass  on $I_{T_j}$ for every $1 \le j \le m$.  The first desired assertion follows from the following claim. \\

\noindent
\textbf{Claim.} Let  $P_j$ be one of currents $T_{1}, \ldots, T_{m}$ for $1 \le j \le m$. Then, the currents $P_1, \ldots, P_{m}$ are of full mass intersection relative to $T$.  
\\

\noindent
Let $\tilde{m}$ be an integer such that  there are at least $\tilde{m}$ currents among $P_{1}, \ldots, P_{m}$ which are equal. We have $0 \le \tilde{m} \le m$. We will prove Claim  by  induction on $\tilde{m}$. When $\tilde{m}=m$, the desired assertion is clear by the hypothesis. Assume that it holds for every $\tilde{m}'> \tilde{m}$. We need to prove it for $\tilde{m}$. By Lemma \ref{le-hasnomassonITjkahler}, we can assume that $\theta_j$'s are all equal to a form $\theta$. 

Without loss of generality, we can assume that $P_{j}=T_1$ for every $1 \le j  \le \tilde{m}$ and $P_{(\tilde{m}+1)}= T_2$. 
If $P_j= T_l$, then we define $v_{jk}:= u_{lk}$ and $P_{jk}:= T_{lk}$.   Put 
$$Q:= \wedge_{j=\tilde{m}+2}^{m} P_{jk} \wedge T\wedge \omega^{n-m-p}, \quad \tilde{P}_{k}:= \ddc \max\{\max\{u_{1}, u_{2}\},-k\}+ \theta$$
and $\tilde{P}:= \ddc \max\{u_{1}, u_{2}\}+ \theta$.  Since the cohomology class of $T_j$ is K\"ahler for every $j$,  we can apply (\ref{ine-dkfullmasscutoffk}) to $P_j$. Hence, in order to obtain the desired assertion,  we need to check that 
\begin{align}\label{ine-T1kmumhaiuphayT2knonweigh}
\int_{B_{k}}T_{1k}^{\tilde{m}} \wedge T_{2k} \wedge Q \to 0
\end{align}
as $k,l \to \infty$, where $$B_{k}:= \{u_{1k} \le -k\}\cup \{ u_{2k} \le -k\}\cup \cup_{j=\tilde{m}+2}^m \{v_{jk} \le -k\}.$$
 By  Theorem \ref{th-equalityMAonstrongerplurifinetopo}, we get
$$\bold{1}_{\{u_{1k} < u_{2k}\}} T_{1k}^{\tilde{m}} \wedge T_{2k} \wedge Q=\bold{1}_{\{u_{1k} < u_{2k}\}} T_{1k}^{\tilde{m}} \wedge \tilde{P}_k \wedge Q.$$
This implies
$$\bold{1}_{\{u_{1k} < u_{2k}\} \cap B_{k}} T_{1k}^{\tilde{m}} \wedge T_{2k} \wedge Q=\bold{1}_{\{u_{1k} < u_{2k}\} \cap B_{k}} T_{1k}^{\tilde{m}} \wedge \tilde{P}_k \wedge Q.$$
It follows that
\begin{align} \label{ine-uocluongtrenu1ku2k}
\int_{ \{u_{1k} < u_{2k}\} \cap B_{k}}T_{1k}^{\tilde{m}} \wedge T_{2k} \wedge Q =\int_{ \{u_{1k} < u_{2k}\} \cap B_{k}}  T_{1k}^{\tilde{m}} \wedge \tilde{P}_k \wedge Q
\end{align}
On the other hand, by  induction hypothesis, the currents $T_1,\ldots T_1, P_{\tilde{m}+2}, \ldots, P_{m}$ ($(\tilde{m}+1)$ times $T_1$) are of full mass intersection relative  to $T$. This combined with Theorem \ref{th-monoticity} implies that  $T_1,\ldots T_1, \tilde{P}, P_{\tilde{m}+2}, \ldots, P_{m}$ ($\tilde{m}$ times $T_1$) are of full mass intersection relative to $T$ because $\tilde{P}$ is less singular than $T_1$. Using this,  (\ref{ine-dkfullmasscutoffk}) and the fact that
$$\{u_{1k} < u_{2k}\} \cap B_{k} \subset   \{u_{1} \le -k\}\cup \{\max\{u_{1},u_{2}\} \le -k\} \cup \cup_{j=\tilde{m}+2}^m \{v_{jk} \le -k\},$$
we see that the right-hand side of (\ref{ine-uocluongtrenu1ku2k}) converges to $0$ as $k \to \infty$.  
It follows that 
\begin{align}\label{ine-chanmhaiphaycong1noweight}
\int_{ \{u_{1k} < u_{2k}\} \cap B_{k}} T_{1k}^{\tilde{m}} \wedge T_{2k} \wedge Q \to 0
\end{align}
as $k \to \infty$.  

Now let $\tilde{P}'_k:= \ddc \max\{u_{1k}, u_{2k}-1\}+ \theta$. The last current converges to 
$$\tilde{P}':= \ddc \max\{u_1, u_2-1\}+ \theta.$$
 Observe that $\tilde{P}'$ is less singular than $T_2$.  By induction hypothesis and Theorem \ref{th-monoticity} as above, we see that the currents $\tilde{P}', \ldots, \tilde{P}', T_2, P_{\tilde{m}+2}, \ldots, P_{m}$ ($\tilde{m}$ times $\tilde{P}'$) are of full mass intersection relative to $T$. Moreover, since $T$ has no mass on $I_{T_1}$, by  $(iii)$ of Proposition \ref{pro-sublinearnonpluripolar}, we get
\begin{align}\label{eq-massIT1Pnga'cnvex}
\bold{1}_{I_{T_1}} \langle \tilde{P}'^m \wedge T_2 \wedge \wedge_{j=\tilde{m}+2}^m P_j \dot{\wedge} T \rangle =0
\end{align}
  Using similar arguments as in the first part of the proof, we obtain
\begin{align} \label{ine-tinhu1lonhonu2}
\int_{ \{u_{1k} > u_{2k} -1\}\cap B_{k}}  T_{1k}^{\tilde{m}} \wedge T_{2k} \wedge Q  &= \int_{ \{u_{1k} > u_{2k} -1\}\cap B_{k}}   \tilde{P}'^{\tilde{m}}_k \wedge T_{2k} \wedge Q \le  \int_{B_{k}}  \tilde{P}'^{\tilde{m}}_k \wedge T_{2k} \wedge Q\\
\nonumber
& = \int_{X}  \tilde{P}'^{\tilde{m}}_k \wedge T_{2k} \wedge Q- \int_{X \backslash B_k}\tilde{P}'^{\tilde{m}}_k \wedge T_{2k} \wedge Q\\
\nonumber
&= \int_{X}  \tilde{P}'^{\tilde{m}}_k \wedge T_{2k} \wedge Q- \int_{X \backslash B_k}\langle \tilde{P}'^m \wedge T_2 \wedge \wedge_{j=\tilde{m}+2}^m P_j \dot{\wedge} T \rangle\\
\nonumber
&= \big\| \langle \tilde{P}'^m \wedge T_2 \wedge \wedge_{j=\tilde{m}+2}^m P_j \dot{\wedge} T \rangle \big\|_{B_{k}}
\end{align}
because $\tilde{P}', \ldots, \tilde{P}', T_2, P_{\tilde{m}+2}, \ldots, P_{m}$ ($\tilde{m}$ times $\tilde{P}'$) are of full mass intersection relative to $T$. Observe that the right-hand side of (\ref{ine-tinhu1lonhonu2})  converges to 
$$\big\| \langle \tilde{P}'^m \wedge T_2 \wedge \wedge_{j=\tilde{m}+2}^m P_j \dot{\wedge} T \big\|_{I_{T_1}}$$
as $k \to \infty$. The last quantity is equal to $0$ because of (\ref{eq-massIT1Pnga'cnvex}). Consequently, we get
\begin{align*}
\int_{ \{u_{1k} > u_{2k}-1\} \cap B_{k}} T_{1k}^{\tilde{m}} \wedge T_{2k} \wedge Q \to 0
\end{align*}
as $k \to \infty$. Combining this and (\ref{ine-chanmhaiphaycong1noweight}) and the fact that  $X= \{u_{1k} < u_{2k}\} \cup \{u_{1k} > u_{2k}-1\}$   gives (\ref{ine-T1kmumhaiuphayT2knonweigh}).  The proof is finished.
\endproof

Let $\chi: \R \to \R$ be a continuous increasing function such that $\chi(-\infty)= -\infty$  and  for every constant $c$,  there exist constants $M,c_1, c_2>0$ so that  
\begin{align}\label{ine-prochiweighted}
\chi(t+c) \ge c_1 \chi(t) - c_2
\end{align}
for $t< -M$. Such a function is called \emph{a weight}.   For every function $\xi$ bounded from above on $X$, let
\begin{align} \label{ine-def-weightedclass}
\tilde{E}_\xi(T_1,\ldots, T_m,T):=\int_X - \xi \langle T_1 \wedge  \cdots \wedge T_m   \dot{\wedge}  T\rangle \wedge \omega^{n-m-p}
\end{align}
The following simple lemma, which generalizes \cite[Proposition 1.4]{GZ-weighted},  will be useful in practice. 

\begin{lemma}\label{le-anotherequivcond-weight} Let  $\xi:= \chi(\sum_{j=1}^m u_j)$ and $\xi_k:= \chi(\max\{\sum_{j=1}^m  u_{j}, -k\})$. Assume that  the current $T_1, \ldots, T_m$ are of full mass intersection relative to $T$. Then we have
\begin{align}\label{eq-tinhzikTjkT2}
\tilde{E}_{\xi_k}(T_{1k}, \ldots, T_{mk}, T)= \int_X - \xi_k \langle \wedge_{j=1}^m T_j  \dot{\wedge}  T \rangle \wedge \omega^{n-m-p}
\end{align}
which increases to $\tilde{E}_\xi(T_1, \ldots,T_m,T)$ as $k \to \infty$; and if additionally $\tilde{E}_\xi(T_1, \ldots,T_m,T)< \infty$, then
\begin{align}\label{eq-tinhzikTjkT}
\big\| \xi_k \wedge_{j=1}^m T_{jk} \wedge T- \xi_k \langle \wedge_{j=1}^m T_j  \dot{\wedge}  T \rangle \big\| \to 0
\end{align} 
as $k \to \infty$.
\end{lemma}

\proof
Put $Q:= \langle \wedge_{j=1}^m T_j  \dot{\wedge}  T \rangle$ and $Q_k:=  \wedge_{j=1}^m T_{jk} \wedge T$. Let $\xi^0_k$ be the value of $\xi_k$ on  $\cup_{j=1}^m  \{u_j \le -k\}$. We have 
\begin{align*}
\xi_k Q_k &=\xi_k \bold{1}_{\cap_{j=1}^m  \{u_j >-k\}} Q + \xi^0_k\bold{1}_{\cup_{j=1}^m  \{u_j \le -k\}} Q_k\\
&=\xi_k \bold{1}_{\cap_{j=1}^m  \{u_j >-k\}} Q + \xi^0_k Q_k - \xi^0_k\bold{1}_{\cap_{j=1}^m  \{u_j > -k\}} Q\\
&=\xi_k Q + \xi^0_k \big(Q_k - Q\big).
\end{align*}
Hence (\ref{eq-tinhzikTjkT2}) follows by integrating the last equality  over $X$ and using the hypothesis. We check (\ref{eq-tinhzikTjkT}). Since
$Q_k=Q$ on $\cap_{j=1}^m  \{u_j >-k\}$, we have
$$ \xi_k Q_k- \xi_k Q=  \bold{1}_{\cup_{j=1}^m  \{u_j \le  -k\}}\xi_k Q_k- \bold{1}_{\cup_{j=1}^m  \{u_j \le -k\}}\xi_k Q.$$
On the other hand, observe that   
$$\|\bold{1}_{\cap_{j=1}^m  \{u_j >-k\}} \xi_k Q \| \to \| \xi Q\|= E_\xi(T_1, \ldots, T_m,T)<\infty$$
 as $k \to \infty$ because $Q$ has no mass on $\cup_{j=1}^m \{u_j = -\infty\}$. Hence $\|\bold{1}_{\cup_{j=1}^m  \{u_j \le -k\}} \xi_k Q \| \to 0$ as $k \to \infty$. Observe
\begin{align*} 
\int_X  \bold{1}_{\cup_{j=1}^m  \{u_j \le -k\}} \xi_k Q_k  \wedge \omega^{n-m-p} &= \int_X  (\xi_k Q_k -  \bold{1}_{\cap_{j=1}^m  \{u_j > -k\}} \xi_k Q_k) \wedge \omega^{n-m-p}\\
&= \int_X  (\xi_k Q_k -  \bold{1}_{\cap_{j=1}^m  \{u_j > -k\}} \xi_k Q) \wedge \omega^{n-m-p}
\end{align*}
converging to $0$ as $k \to \infty$ by (\ref{eq-tinhzikTjkT2}). Thus (\ref{eq-tinhzikTjkT}) follows.
 The proof is finished.
\endproof

\begin{definition} \label{def-weigtedclasschi} We say that $T_1, \ldots, T_m$ are \emph{of full mass intersection relative to $T$ with weight $\chi$} if  $T_1, \ldots, T_m$ are of full mass intersection relative to $T$ and  for every nonempty set $J \subset \{1, \ldots, m\}$, we have $\tilde{E}_\xi\big((T_j)_{j \in J},T\big)< \infty$, where $\xi:= \chi(\sum_{j=1}^m u_j)$ .
\end{definition}

The last definition is independent of the choice of potentials $u_j$ by (\ref{ine-prochiweighted}).   For  currents $T_1, \ldots, T_m$ of full mass intersection relative to $T$ and $\xi:= \chi(\sum_{j=1}^m u_j)$, we can also define \emph{the joint $\xi$-energy relative to $T$} of $T_1, \ldots, T_m$ by putting 
$$E_\xi(T_1, \ldots, T_m,T):= \sum_{J} \int_X - \xi \langle \wedge_{j \in J} T_j  \dot{\wedge}  T \rangle \wedge \omega^{n-p- |J|},$$ 
where the sum is taken over every subset $J$ of $\{1, \ldots, m\}$. The last energy depends on the choice of potentials but its finiteness does not.

Let $\alpha_1, \ldots, \alpha_m$ be K\"ahler $(1,1)$-classes on $X$. Let $\mathcal{E}_{\chi}(\alpha_1, \ldots, \alpha_m, T)$ be  the set of  $m$-tuple $(T_1, \ldots,T_m)$ of closed positive $(1,1)$-currents such that $T_1, \ldots, T_m$ of full mass intersection relative to $T$ with weight $\chi$ and  $T_j \in \alpha_j$ for $1 \le j \le m$. 

For pseudoeffective $(1,1)$-class $\beta$, we define $\mathcal{E}_{\chi, m}(\beta, T)$ to be the subset of  $\mathcal{E}_m (\beta, T)$ consisting of $P$ such that $(P, \ldots, P)$ is in $\mathcal{E}_\chi(\beta, \ldots, \beta, T)$  ($m$ times $\beta$).  When $m=n$, the  class $\mathcal{E}_{\chi,m}(\beta, T)$ was mentioned in \cite[Section 5.2.2]{GZ-weighted}. The last class was studied in \cite{GZ-weighted,BEGZ} when $T$ is the current of integration along $X$ and $m=n$. 

For $P \in  \mathcal{E}_{\chi, m}(\beta, T)$ with $P= \ddc u+ \theta$ ($\theta$ is K\"ahler) and $\xi:= \chi(m\, u)$, we put 
$$E_{\xi}(P,T):= E_\xi (P, \ldots, P, T)$$
 ($m$ times $P$). The following result explains why our weighted class generalizes that given in \cite{GZ-weighted}. 

\begin{lemma} \label{le-generalwightedclassGZ} $(i)$ A current   $P \in \mathcal{E}_{\chi, m}(\beta, T)$ if and only if  $P \in \mathcal{E}_m(\beta,T)$ and $\xi$ is integrable with respect to $\langle P^m \wedge T \rangle$.

$(ii)$ Assume that $T$ is the current of integration along $X$ and $m=n$. Then,  the energy $E_{\chi, m}(P,T)$ is equivalent to the energy associated to $\tilde{\chi}(t):= \chi(n \, t)$ given in \cite{GZ-weighted}.
\end{lemma}

\proof This is a direct consequence of  computations in the proof of \cite[Proposition 2.8 (i)]{BEGZ} and Lemma \ref{le-anotherequivcond-weight}: one first consider the case where   $\chi$ is smooth and use Lemma \ref{le-tinhddcchivvoichiC2} below;  the general case follows by regularizing $\chi$. The proof is finished. 
\endproof

 The following result is obvious. 

\begin{lemma}\label{le-congthemthetajchobang} Let $\eta_j$ be a positive closed $(1,1)$-form for $1 \le j \le m$. Then  $T_1, \ldots, T_m$ are of full mass intersection relative to $T$ with weight $\chi$ if and only if $(T_1+ \eta_1), \ldots, (T_m+\eta_m),T$ are of full mass intersection relative to $T$ with weight $\chi$.
\end{lemma}

From now on, we focus on convex weights.  Let $\mathcal{W}^-$ be the set of convex increasing functions $\chi: \R \to \R$ with $\chi(-\infty)= - \infty$.  Observe that if $\chi \in \mathcal{W}^- $, then $\chi$ is automatically continuous.    Basic examples of $\chi$ are $-(-t)^r$ for $0 < r \le 1$. We have the following important observation which, in particular, implies that $\mathcal{W}^-$ is a set of weights in the sense given above.

\begin{lemma} \label{le-regularizedchi} Let $\chi \in \mathcal{W}^-$.  Let $g$ be a smooth radial cut-off function on $\R$, \emph{i.e,} $g(t)= g(-t)$ for $t \in \R$, $g$ is of compact support ,  $0 \le g \le 1$ and $\int_\R g(t) dt =1$. Put $g_\epsilon(t):=  \epsilon^{-1}g(\epsilon t)$ for every constant $\epsilon >0$ and $\chi_\epsilon:= \chi * g_\epsilon$ (the convolution of $\chi$ with $g_\epsilon$). Then   $\chi_\epsilon \in \mathcal{W}^-$, $\chi_\epsilon \ge \chi$ and 
\begin{align}\label{def-Wtru}
0 \ge t \chi'_\epsilon(t) \ge \chi_\epsilon(t) - \chi_\epsilon(0)
\end{align}
 for $t\le 0$.  Consequently,  (\ref{ine-prochiweighted}) holds for $\chi \in \mathcal{W}^- $, or in other words, $\chi$ is a weight. 
\end{lemma}

Clearly we always have  that $\chi_\epsilon$ converges uniformly to $\chi$ because of continuity of $\chi$. 

\proof By definition, 
$$\chi_\epsilon(t)= \int_\R \chi(t-s) g_\epsilon(s) ds= \int_{\R^+}\big( \chi(t-s)+ \chi(t+s)\big) g_\epsilon(s) ds.$$ 
We deduce that since $\chi$ is convex, $\chi_\epsilon \ge \chi$ for every $\epsilon$. The inequality (\ref{def-Wtru}) is a direct consequence of the convexity.     It remains to prove the last desired assertion. We can assume $\chi$ is smooth by the previous part of the proof. Since $\chi$ is increasing, it is enough to prove the desired assertion for $c\ge 0$.  Fix a constant $c\ge 0$. Consider $t < -|c|+M$ for some big constant $M$. Write $\chi(t+c)= \chi(t)+ \int_t^{t+c} \chi'(r) dr$. Using (\ref{def-Wtru}), we obtain that  there is a constant $c_1>0$ independent of $c$ such that
$$\chi(t+c)- \chi(t) \ge  - \int_{|t|}^{|t|+ c} [-\chi(-r)+ c_1]/r dr \ge   ( \chi(t)-c_1)\big(\log(|t|+c)- \log |t|\big).$$
Thus the desired assertion follows.  The proof is finished.
\endproof

We will need the following computation which seems to be used implicitly in the literature. 

\begin{lemma} \label{le-tinhddcchivvoichiC2}  Let $\chi \in \cali{C}^3(\R)$ and $w_1, w_2$  bounded psh functions on an open subset $U$ of $\C^n$. Let $Q$ be a closed positive current of bi-dimension $(1,1)$ on $U$. Then we have 
\begin{align}\label{eq-ddcchi(w2)integbyparta}
\ddc \chi(w_2) \wedge Q = \chi''(w_2) d w_2 \wedge \dc w_2 \wedge Q+ \chi'(w_2) \ddc w_2 \wedge Q
\end{align}
and the operator $w_1 \ddc \chi(w_2) \wedge Q$ is continuous (in the usual weak topology of currents) under decreasing sequences of smooth psh functions converging to $w_1,w_2$.
Consequently, if $f$ is a smooth function with compact support in $U$, then the equality
\begin{align}\label{eq-ddcchi(w2)integbyparta2}
\int_U f w_1 \ddc \chi(w_2) \wedge Q= \int_U \chi(w_2) \ddc (f w_1) \wedge Q
\end{align}
holds. Moreover, for $f$ as above, we also have 
\begin{align}\label{eq-ddcchi(w2)integbyparta23}
\int_U f  \chi(w_2) \ddc w_1  \wedge Q= \int_U \chi(w_2) df \wedge  \dc w_1 \wedge Q+ \int_U f \chi'(w_2) d w_2 \wedge  \dc w_1 \wedge Q.
\end{align}
\end{lemma}

\proof Clearly, all of three desired equalities follows from the integration by parts if $w_1, w_2$ are smooth.  The arguments below essentially say that both sides of these equalities are continuous under sequences of smooth psh functions decreasing to $w_1, w_2$. This is slightly non-standard due to the presence of $Q$ even when $\chi$ is convex. 

First observe that (\ref{eq-ddcchi(w2)integbyparta2}) is a consequence of the second desired assertion because both sides of (\ref{eq-ddcchi(w2)integbyparta2}) are continuous under a sequence of smooth psh functions decreasing to $w_2$. We prove (\ref{eq-ddcchi(w2)integbyparta}).   The desired equality (\ref{eq-ddcchi(w2)integbyparta}) clearly holds if $w_2$ is smooth.  In general, let $(w_2^\epsilon)_\epsilon$ be a sequence of standard regularisations of $w_2$.  
Recall that $ \ddc \chi(w_2) \wedge Q $ is defined to be $\ddc \big( \chi(w_2)Q\big)$ which is equal to the limit of $\ddc \big( \chi(w^\epsilon_2)Q\big)$ as $\epsilon \to 0$.  By (\ref{eq-ddcchi(w2)integbyparta}) for $w_2^\epsilon$ in place of $w_2$, we see that $\ddc \big( \chi(w^\epsilon_2)Q\big)$ is of uniformly bounded mass. As a result, $\ddc \chi(w_2) \wedge Q$ is of order $0$. Thus $w_1 \ddc \chi(w_2)\wedge Q$ is well-defined. Put  
$$I(w_1, w,w_2):= w_1 \chi''(w) d w_2 \wedge \dc w_2 \wedge Q+ w_1\chi'(w) \ddc w_2 \wedge Q.$$
Recall that $I(1, w_2^\epsilon, w_2^\epsilon) \to \ddc \chi(w_2) \wedge Q$.   By Corollary \ref{cor-convergplurifine}, we have
\begin{align}\label{eq-limitw2epsilon}
I(w_1,w_2, w_2^\epsilon) \to I(w_1,w_2,w_2)
\end{align}
as $\epsilon \to 0$. On the other hand, since $\chi''$ is in $\cali{C}^1$, we get 
$$|\chi''(w_2^\epsilon)- \chi''(w_2) | \lesssim (w_2^\epsilon- w_2), \quad |\chi'(w_2^\epsilon)- \chi'(w_2) | \lesssim (w_2^\epsilon- w_2).$$
 This combined with the convergence of Monge-Amp\`ere operators under decreasing sequences tells us that
\begin{align}\label{eq-limitw2epsilon2}
\big(I(w_1,w_2^\epsilon, w_2^\epsilon)- I(w_1,w_2, w_2^\epsilon) \big) \to 0
\end{align}
as $\epsilon \to 0$. Combining (\ref{eq-limitw2epsilon2}) and (\ref{eq-limitw2epsilon}) gives that $I(w_1,w_2^\epsilon, w_2^\epsilon) \to I(w_1,w_2,w_2)$ as $\epsilon \to 0$. Letting $w_1\equiv 1$ in the last limit, we get (\ref{eq-ddcchi(w2)integbyparta}). The second desired assertion also follows.  We prove (\ref{eq-ddcchi(w2)integbyparta23}) similarly.  The proof is finished.  
\endproof


Here is a  monotonicity property of weighted classes.

\begin{theorem} \label{the-monotonicityforweigtedclass} (Monotonicity of weighted classes)  Let $\chi \in \mathcal{W}^-$ with $|\chi(0)| \le 1$.  Let  $T'_j $ be  a closed positive $(1,1)$-current whose cohomology class is K\"ahler  for $1 \le j \le m$. Assume that $T'_j$ is less singular than $T_j$ and $T_1, \ldots, T_m$ are of full mass intersection  relative to $T$  with weight $\chi$. Then, $T'_1, \ldots, T'_m$ are also of full mass intersection relative to $T$ with weight $\chi$ and  for $\xi:=  \chi\big(\sum_{j=1}^m u_j\big)$, we have   
\begin{align}\label{ine-ExiTphayjmonotweightedclass}
E_{\xi}(T'_1, \ldots, T'_m,T) \le c_1 E_\xi(T_1, \ldots, T_m,T)+c_2,
\end{align}
for some constants $c_1,c_2>0$ independent of $\chi$. 
\end{theorem}


\proof   Let $T'_j= \ddc u'_j+ \theta'_j$ for some smooth K\"ahler forms  $\theta'_j$. Put  $v:= \sum_{j=1}^m u_j$, $v_k:= \max\{v, -k\}$,   $\xi:= \chi(v)$ and  $\xi_k:= \chi(v_k)$. Let $u_{jk}:= \max\{u_j, -k\}$ ,  $T_{jk}:= \ddc u_{jk}+ \theta_j$.  Define $u'_{jk}, T'_{jk}$ similarly. We have $u_j \le u'_j$.  Observe that
\begin{align}\label{inclu-weightedTphayTtrongfulllmass}
(T'_j)_{j \in J'}, (T_j)_{j \in J} \,  \,\text{are of full mass intersection relative to $T$}
\end{align}
 by Lemma \ref{le-hasnomassonITjkahler}. Hence, in order to obtain the first desired assertion, it suffices to check (\ref{ine-ExiTphayjmonotweightedclass}) because $- \chi(\sum_{j=1}^m u'_j) \le - \xi$.  The desired inequality  (\ref{ine-ExiTphayjmonotweightedclass})  follows by letting $k \to \infty$ in  the following claim and using  Lemma \ref{le-anotherequivcond-weight} and (\ref{inclu-weightedTphayTtrongfulllmass}).\\

\noindent
\textbf{Claim.}  For every $J, J' \subset \{1, \ldots, m\}$ with $J \cap J' = \varnothing$, we have  
$$\int_X \xi_k   \wedge_{j \in J'} T'_{jk} \wedge_{j \in J} T_{jk} \wedge  T \wedge \omega^{n-p-|J|- |J'|}> - c_1 E_{\xi_k}(T_{1k}, \ldots, T_{mk},T)- c_2,$$
for some constants  $c_1, c_2$ independent of $\chi$. \\

\noindent
It remains to prove Claim now.  We observe  that it is enough to prove Claim for $\chi$ smooth by Lemma \ref{le-regularizedchi}. Consider, from now on, \emph{smooth} $\chi$. We can also assume that $u_j,u'_j<-M$ for some big constant $M$.   

We prove Claim by induction on $|J'|$.  When $J'= \varnothing$, this is clear. Assume Claim holds for every $J'$ with $|J'|< m'$. We need to prove that it holds for $J'$ with $|J'|=m'$. Without loss of generality, we can assume that $ 1 \in J'$. 
 Put
$$Q_k:= \wedge_{j \in J' \backslash \{1\}} T'_{jk} \wedge \wedge_{j \in J}T_{jk} \wedge T \wedge \omega^{n- p- |J|- |J'|}, \quad Q:= \langle \wedge_{j \in J' \backslash \{1\}} T'_{j} \wedge_{j \in J}T_{j} \dot{\wedge} T \rangle \wedge \omega^{n- p- |J|- |J'|}.$$
  By integration by parts, we obtain 
\begin{align} \label{eq-xikT1kQkeweighted}
I_k:=  \int_X \xi_k  T'_{1k}\wedge Q_k=  \int_X u'_{1k} \ddc \xi_k  \wedge Q_k+  \int_X \xi_k \theta'_1 \wedge Q_k.
\end{align}
Denote by $I_{k,1}, I_{k,2}$ the first and second terms in the right-hand side of the last equality. By induction hypothesis, we get
$$I_{k,2} \ge - c_1 E_{\xi_k}(T_{1k}, \ldots, T_{mk},T)- c_2 $$
for some constants $c_1, c_2>0$ depending only on the masses of $T_j,T'_j$ for $1 \le j \le m$. 

It remains to treat $I_{k,1}$.   Put $\theta:= \sum_{j=1}^m \theta_j$. Note that $\ddc v_k + \theta \ge 0$.   By Lemma \ref{le-tinhddcchivvoichiC2}, the current $(\ddc \xi_k+ \chi'(v_k) \theta)  \wedge Q_k$ is positive. This combined with the inequality $u_{1k} \le u'_{1k}$  gives
\begin{align*}
I_{k,1}  & =  \int_X u'_{1k} (\ddc  \xi_k+ \chi'(v_k) \theta)  \wedge Q_k -  \int_X u'_{1k}\chi'(v_k)  \theta \wedge Q_k \\
&\ge  \int_X u_{1k} (\ddc  \xi_k+ \chi'(v_k) \theta)  \wedge Q_k -  \int_X u'_{1k}\chi'(v_k)  \theta  \wedge Q_k \\
& =  \int_X   \xi_k T_{1k} \wedge Q_k- \int_X  \xi_k \theta_1 \wedge Q_k  +  \int_X (u_{1k}- u'_{1k}) \chi'(v_k) \theta \wedge Q_k\\
&\ge \int_X  \xi_k T_{1k} \wedge Q_k +\int_X  \xi_k (\theta-\theta_1) \wedge Q_k - \chi(0) \int_X \theta \wedge Q_k
\end{align*}
because  
$$(u_{1k}- u'_{1k}) \chi'(v_k)  \ge u_{1k} \chi'(v_k)\ge v_k \chi'(v_k) \ge  \chi(v_k) - \chi(0)= \xi_k - \chi(0)$$
 (see (\ref{def-Wtru})). This combined with induction hypothesis gives
$$I_{k,2} \ge - c_1 E_{\xi_k}(T_{1k}, \ldots, T_{mk},T)- c_2.$$
Hence, Claim  follows.  The proof is finished.
\endproof


By the above proof, one can check that if we fix $\theta_1, \ldots, \theta_m$ and $\theta'_1, \ldots, \theta'_m$, then the constants $c_1$ and $c_2$ in (\ref{ine-ExiTphayjmonotweightedclass}) can be chosen to be independent of $T_1, \ldots, T_m$.  Here is a convexity property for weighted classes. 

\begin{theorem}  \label{the-mathcalEnalphaweighted} Let $\chi \in \mathcal{W}^-$. Assume that  $T_j, \ldots, T_j$  ($m$ times $T_j$) are of full mass intersection relative to $T$ with weight $\chi$   for $1 \le j \le m$. Then,  $T_1, \ldots, T_m$ are also of full mass intersection relative to $T$ with weight $\chi$. In particular,  for  K\"ahler  $(1,1)$-classes $\beta, \beta'$, we have
 $$\mathcal{E}_{\chi, m}(\beta, T)+ \mathcal{E}_{\chi, m}(\beta',T) \subset \mathcal{E}_{\chi,m}(\beta+ \beta', T)$$
and $\mathcal{E}_{\chi, m}(\beta,T)$ is convex. 
\end{theorem}

\proof  The second and third desired assertions are  direct consequences of  the first one. We prove the first desired assertion. Firstly, by Theorem \ref{the-mathcalEnalpha}, we have 
\begin{align}\label{inclu-P1PsT}
T_1, \ldots, T_m \, \, \text{are of full mass intersection relative to $T$}.
\end{align} 
By Lemma \ref{le-congthemthetajchobang}, we can assume that $\theta_j= \theta$ for every $1 \le j \le m$. The desired assertion is a consequence of the following claim.\\

\noindent
\textbf{Claim.} Let  $P_j$ be one of currents $T_{1}, \ldots, T_{m}$ for $1 \le j \le m$. Then, the currents $P_1, \ldots, P_{m}$ are of full mass intersection relative to $T$ with weight $\chi$.  
\\

\noindent
Let $\tilde{m}$ be an integer such that  there are at least $\tilde{m}$ currents among $P_{1}, \ldots, P_{m}$ which are equal. We have $0 \le \tilde{m} \le m$. We will prove Claim  by  induction on $\tilde{m}$. When $\tilde{m}=m$, the desired assertion is clear by the hypothesis. Assume that it holds for every number $\tilde{m}'> \tilde{m}$. We need to prove it for $\tilde{m}$. 

Without loss of generality, we can assume that $P_{j}=T_1$ for every $1 \le j  \le \tilde{m}$ and $P_{(\tilde{m}+1)}= T_2$. 
If $P_j= T_l$, then we define $v_{jk}:= u_{lk}$ and $P_{jk}:= T_{lk}$.   Put $Q:= \wedge_{j=\tilde{m}+2}^{m} P_{jk} \wedge T\wedge \omega^{n-m-p}$,  $\tilde{P}_{k}:= \ddc \max\{u_{1k}, u_{2k}\}+ \theta$ and $$w_k:= \sum_{j=\tilde{m}+2}^{m} v_{jk}.$$
By Lemma \ref{le-anotherequivcond-weight} and (\ref{inclu-P1PsT}),  we need to check that 
\begin{align}\label{ine-T1kmumhaiuphayT2k}
\int_X - \chi(\tilde{m} u_{1k}+ u_{2k}+ w_k) T_{1k}^{\tilde{m}} \wedge T_{2k} \wedge Q \le C
\end{align}
for some constant $C$ independent of $k$. Actually, we need to verify a stronger statement that for every $J \subset \{1,\ldots, \tilde{m}\}$, then $$\int_X - \chi(\tilde{m} u_{1k}+ u_{2k}+ w_k) \wedge_{j \in J}P_{jk} \wedge T \wedge \omega^{n-p-|J|}$$
is uniformly bounded in $k$, but the proof will be similar to that of (\ref{ine-T1kmumhaiuphayT2k}). 

Observe that $X= \{u_{1k} < u_{2k}\} \cup \{u_{1k} > u_{2k}-1\}$.  Using Theorem \ref{th-equalityMAonstrongerplurifinetopo}, we get
$$\bold{1}_{\{u_{1k} < u_{2k}\}} T_{1k}^{\tilde{m}} \wedge T_{2k} \wedge Q=\bold{1}_{\{u_{1k} < u_{2k}\}} T_{1k}^{\tilde{m}} \wedge \tilde{P}_k \wedge Q.$$
This combined with the fact that $ - \chi(\tilde{m} u_{1k} + u_{2k}+ w_k) \le - \chi\big((\tilde{m}+1) u_{1k}+ w_k\big)$ on  $\{u_{1k} < u_{2k}\}$ implies
\begin{multline*}
\int_{ \{u_{1k} < u_{2k}\}} - \chi(\tilde{m} u_{1k} + u_{2k}+ w_k) T_{1k}^{\tilde{m}} \wedge T_{2k} \wedge Q =\\
\int_{ \{u_{1k} < u_{2k}\}} - \chi(\tilde{m} u_{1k} + u_{2k}+ w_k) T_{1k}^{\tilde{m}} \wedge \tilde{P}_k \wedge Q
\end{multline*}
which is $\le$
 $$ \int_{X} -  \chi\big((\tilde{m}+1) u_{1k} + w_k\big) T_{1k}^{\tilde{m}} \wedge \tilde{P}_k \wedge Q.$$
By the fact that $\tilde{P}_k$ is less singular than $T_{1k}$ and  Claim in the proof of Theorem \ref{the-monotonicityforweigtedclass}, we see that the right-hand side of the last inequality is bounded uniformly in $k$ because $T_1,\ldots T_1, P_{\tilde{m}+2}, \ldots, P_{m}$ ($(\tilde{m}+1)$ times $T_1$) are of full mass intersection relative to $T$ with weight $\chi$ (this is a consequence of the induction hypothesis). It follows that 
\begin{align}\label{ine-chanmhaiphaycong1}
\int_{ \{u_{1k} < u_{2k}\}} - \chi(\tilde{m} u_{1k} + u_{2k}+ w_k) T_{1k}^{\tilde{m}} \wedge T_{2k} \wedge Q \le C
\end{align}
for some constant $C$ independent of $k$. Similarly, for $\tilde{P}'_k:= \ddc \max\{u_{1k}, u_{2k}-1\}+ \theta$,  we have 
\begin{multline*}
\int_{ \{u_{1k} > u_{2k} -1\}} - \chi(\tilde{m} u_{1k} + u_{2k}+ w_k) T_{1k}^{\tilde{m}} \wedge T_{2k} \wedge Q  \le \\
 c_1\int_{X} -  \chi\big((\tilde{m}+1) u_{2k} + w_k\big)(\tilde{P}_k')^{\tilde{m}} \wedge T_{2k} \wedge Q +c_2
\end{multline*}
for some constants $c_1, c_2>0$ independent of $k$ (we used (\ref{ine-prochiweighted}) here). Using the induction hypothesis again and the fact that $\tilde{P}'_k$ is less singular that $T_{2k}$, we obtain 
\begin{align}\label{ine-chanmhaiphaycong2}
\int_{ \{u_{1k} > u_{2k} -1\}} - \chi(\tilde{m} u_{1k} + u_{2k}+ w_k) T_{1k}^{\tilde{m}} \wedge T_{2k} \wedge Q \le C
\end{align}
for some constant $C$ independent of $k$. Combining (\ref{ine-chanmhaiphaycong2}) and (\ref{ine-chanmhaiphaycong1}) gives (\ref{ine-T1kmumhaiuphayT2k}). The proof is finished.
\endproof

We  now give  a continuity property of weighted classes which generalizes the second part of \cite[Theorem 2.17]{BEGZ}) in the case where the cohomology classes of currents are K\"ahler.

\begin{proposition} \label{pro-continuityweightedclass} Let $\chi \in \mathcal{W}^-$.  Let $T_{j,l}= \ddc u_{j,l}+ \theta_j$ be closed positive $(1,1)$-currents, for every $1 \le j \le m$, $l \in \N$ such that $u_{j,l} \to u_j$ as $l \to \infty$ in $L^1$ and  $u_{j,l}\ge u_j$ for every $j,l$. Assume that  $T_1, \ldots, T_m$ are of full mass intersection relative to $T$ with weight $\chi$. Then, we have 
$$ \chi\big(\sum_{j=1}^m u_{j,l}\big)\langle \wedge_{j=1}^m T_{j ,l} \dot{\wedge}  T\rangle \to \chi\big(\sum_{j=1}^m u_{j}\big)\langle\wedge_{j=1}^m T_{j} \dot{\wedge}  T\rangle$$  
as $l \to \infty$.
\end{proposition}

\proof By Theorem \ref{the-monotonicityforweigtedclass}, the currents $T_{1,l}, \ldots, T_{m,l}$ are of full mass intersection relative to $T$ with weight $\chi$ for every $l$. The desired convergence now follows by using arguments similar to those in the proof of \cite[Theorem 2.17]{BEGZ} (notice also Lemma \ref{le-anotherequivcond-weight} and the comment after Theorem \ref{the-monotonicityforweigtedclass}).
The proof is finished.
\endproof

\begin{remark}\label{re-convexity}   As in \cite{GZ-weighted}, we can check that 
$\mathcal{E}(\alpha_1, \ldots, \alpha_m, T)= \cup_{\chi \in \mathcal{W}^-}\mathcal{E}_\chi(\alpha_1, \ldots, \alpha_m, T)$
for K\"ahler classes $\alpha_1, \ldots \alpha_m$ on $X$.
\end{remark}

\bibliography{biblio_family_MA,biblio_Viet_papers}
\bibliographystyle{siam}

\bigskip

\noindent
\Addresses
\end{document}